\documentclass[12pt,a4paper,twoside]{article}
\usepackage[utf8]{inputenc}
\usepackage[T1]{fontenc}
\usepackage[margin=2.5cm]{geometry}
\usepackage{lmodern}
\usepackage[french,english]{babel}
\usepackage{bef_aras_nothm}
\usepackage{amsmath,amsthm,amssymb}
\usepackage{fixmath}
\usepackage{upgreek}
\usepackage{enumitem}
\usepackage{amsfonts}
\usepackage{microtype}
\usepackage{float}
\usepackage{mathrsfs}

\usepackage{hyperref}
\hypersetup{
     colorlinks   = true,
     citecolor    = blue
}
\makeatletter
\newcommand{\oset}[3][0ex]{%
    \mathrel{\mathop{#3}\limits^{
            \vbox to#1{\kern-2\ex@
                \hbox{$\scriptstyle#2$}\vss}}}}
\makeatother

\usepackage{verbatim}

\usepackage[headsepline]{scrlayer-scrpage}

\clearscrheadfoot
\ohead{\pagemark}
\ihead{\headmark}
\automark[subsection]{section}

\theoremstyle{plain}

\allowdisplaybreaks[1]

\numberwithin{equation}{section}

\newtheorem{thm}{Theorem}[section]
\newtheorem{lem}[thm]{Lemma}

\theoremstyle{definition}

\newtheorem{rem}[thm]{Remark}

\begin{document}

\newpage
\setcounter{page}{1}
\clearscrheadfoot
\ohead{\pagemark}
\ihead{\headmark}

\title{Abstract nonlinear evolution inclusions of second order with applications in visco-elasto-plasticity}

\author{Aras Bacho\footnotemark[2]}

\date{}
\maketitle

\footnotetext[1]{Ludwig-Maximilians-Universit\"{a}t M\"{u}nchen, Mathematisches Institut, Akademiestr{\ss}e 7, 80333 M\"{u}
nchen, Germany.}

\begin{abstract}
Existence of strong solutions of an abstract Cauchy problem  for a class of doubly nonlinear evolution inclusion of second order is established via a semi-implicit time discretization method. The principal parts of the operators acting on $u$ and $u'$ are multi-valued subdifferential operators and are discretized implicitly. A non-variational and non-monotone perturbation acting nonlinearly on $u$ and $u'$ is allowed and discretized explicitly in time. The convergence of a variational approximation scheme is established using methods from convex analysis. In addition, it is proven that the solution satisfies an energy-dissipation equality. Applications of the abstract theory to various examples, e.g., a model in visco-elastic-plasticity, are provided.
\end{abstract}
\vspace*{1em} \textbf{Keywords} Evolution inclusion of second order $ \cdot $ Nonlinear damping $ \cdot $ Nonsmooth analysis $ \cdot $ Variational approximation scheme  $ \cdot $ Rate-independent dissipation $ \cdot $ Visco-Elasto-Plasticity $ \cdot $ Martensitic transformation in shape-memory alloys 

\section{Introduction}
\subsection{Problem setting}

In this article, we investigate the abstract \textsc{Cauchy} problem 
\begin{align} 
\label{eq:I.2}
\begin{cases}
u''(t)+\partial\Psi_{u(t)}(u'(t))+\partial \calE_t(u(t))+B(t,u(t),u'(t)) \ni f(t), &\quad \text{for a.e. } t\in (0,T),\\
u(0)=u_0, \quad u'(0)=v_0,
\end{cases}
\end{align} where $\Psi_u$ denotes the dissipation potential, $\calE_t$ the energy functional, $B$ the perturbation, and $f$ the external force.
Here, the dissipation potential $\Psi_u$ is, in general, nonlinear, non-quadratic, nonsmooth, and depends nonlinearly on the state $u$. The energy functional $\calE_t=\calE^1+\calE^2_t$ is the sum of a functional $\calE^1$ that is defined by a strongly positive, symmetric, and bounded bilinear form and a strongly continuous functional $\calE^2_t$ such that $\calE_t$ is $\lambda$-convex. The perturbation $B$ is a strongly continuous perturbation of $\partial \Psi_u$ and $\partial\calE_t$. 
\subsection{Illustrative examples}
In the following, we give some illustrative examples of evolution inclusions that can be solved with our abstract theory.\\\\ 
1. In the first example, we consider a visco-elasto-plastic model for the martensitic phase transformation in shape-memory alloys governed by the following system of equations
\begin{align*}
\rho\partial_{tt}\uu+\nu(-1)^n \Delta^n \partial_t\uu -\nabla\cdot(\ssigma_p+\mathbold{\sigma}(\nabla \uu))+\mu (-1)^m \Delta^m \uu = \mathbold{f},\quad \text{in } \Omega_T,\\
\ssigma_p\in \mathrm{Sgn}\left( \lambda'(\nabla \uu(\xx,t)): \nabla\partial_t\uu(\xx,t) \right) \lambda'(\nabla \uu(\xx,t))\quad \text{a.e. in } \Omega_T,
\end{align*} 
2. In the second example, we consider an evolution inclusion with nonlinear damping given by 
\begin{align*}
\partial_{tt} u-\nabla \cdot \mathbf{p} -\Delta u +b(u) = f \quad \text{in } \Omega_T,\\
\mathbf{p}(\xx,t)\in \partial_v \psi(\xx,u(\xx,t), \nabla \partial_t u(\xx,t)) \quad \text{a.e. in } \Omega_T,\end{align*} that for $\psi=0$ and $b(u)=\gamma u, \gamma>0,$ reduces to the classical \textsc{Klein--Gordon} equation, which is a relativistic wave equation with applications in relativistic quantum mechanics. For $\psi(v)=\frac{1}{q}\vert \nabla v\vert^p$, the inclusion can be interpreted as a viscous regularization of the \textsc{Klein--Gordon} equation. \\\\
3. In the final example, we consider the evolution inclusion
\begin{align*}
\partial_{tt} u+\left \vert \partial_t u\right \vert^{q-2}\partial_t  u +p -\nabla \cdot \left( E\nabla u\right) +W'(u) = f \quad \text{in } \Omega_T,\\
p(\xx,t)\in \mathrm{Sgn}\left(\partial_t u(\xx,t)\right)\quad\text{a.e. in } \Omega_T.
\end{align*} which can be interpreted as a model in ferro-magnetism \cite{MiRoSa13NADN}.

In Section \ref{se:App}, we discuss the preceding examples more in detail and show the existence of weak solutions satisfying an energy-dissipation balance. 

\subsection{Literature review}
\label{lit.review}

Very few authors have studied evolution equations of second order where the operator acting on the time derivative of the solution is nonlinear. \textsc{Lions} and \textsc{Strauss} \cite{LioStr65SNLE} showed in their seminal work the well-posedness of the \textsc{Cauchy} problem for the doubly nonlinear evolution equation 
 \begin{align}\label{eq:Cauchy}
u''(t)+A(t)u'(t)+ B(t) u(t)= f(t), \quad t\in (0,T),
\end{align} where $B$ is an unbounded, self-adjoint, and linear operator and $A$ is a nonlinear operator. Under sufficient regularity conditions on the given data and the time dependence of $A$, the authors show well posedness of the problem for two cases with two different methods: compactness and monotonicity methods. The peculiarity in both cases is the assumption that the operators $A(t)$ and $B(t)$ are, for each $t\in[0,T]$, defined on different spaces whose intersection is densely and continuously embedded in both spaces. This implies that the solution $u$ takes values in a different space than its time derivative $u'$. Based on the techniques used in \cite{LioStr65SNLE}, the authors in \cite{EmSiTh15FDNS} showed the existence of solutions to the Cauchy problem for \eqref{eq:Cauchy},  where for each $t\in [0,T]$, $A(t):V_A\rightarrow V_A^*$ is a hemicontinuous operator that satisfy a suitable growth condition such that $A+\kappa I$ is monotone and coercive, and the operator $B(t)=B_0+C(t):V_B\rightarrow V_B^*$ is the sum of a linear, bounded, symmetric, and strongly positive operator and a strongly continuous perturbation $C(t)$. As in \cite{LioStr65SNLE}, the authors assume neither that $V_A$ is continuously embedded in $V_B$ nor the reverse case. The assumptions on $A$ imply that $A+\kappa I$ is maximal monotone and therefore not necessarily a potential operator. Therefore, the result obtained here only partially generalizes the above mentioned results. However, to the best of the authors' knowledge, results on the existence of strong solutions for multivalued operators $A$ which are  nonlinear in $u$ and $u'$ do not exist in the literature. Evolution inclusions occur in many applications, e.g., physical phenomena where rate-independent responses of the body are typical such as in plasticity \cite{MieRou15RIST}, in ferromagnetic hysteresis \cite{Visi00FHO, MiRoSa13NADN} or in Visco-Elasto-Plasticity \cite{RajRou03EDSA}. Applications are also found in optimal control theory \cite{AubCel84DI} or in nonsmooth dynamical systems  \cite{Kunz00NSDS}. Another motivation for this work is to complement the results obtained in B. \cite{Bach25DNEI}, where the principal part of the operator acting on $u$ is nonlinear and multi-valued and the principal part of the operator acting on $u$ is linear, symmetric and positive. Moreover, \cite{Bach25DNEI} allows the operators acting on $u$ and $u'$ to be multi-valued by employing regularization techniques developed in \cite{Bach23GMYR}. A similar extension could be considered in the present work; however, for simplicity, we do not pursue it here and refer the interested reader to \cite{Bach25DNEI}. Our contributions concern the following:
\begin{itemize}
\item We allow the functionals that are acting on $u$ and $u'$ to be nonsmooth, hence generating a nonlinear  multi-valued subdifferential in the equations. 
\item We allow the multi-valued operators to live on different spaces.
\item We extend the applications of the abstract theory for evolution inclusions to important applications in physics, e.g., in models for visco-elasto-plasticity, see Section \ref{ex.martensitic}.
\item We allow non-variational and non-monotone perturbations of the subdifferential operators. 
\end{itemize}
For further results on nonlinear abstract evolution inclusions, we refer to \cite{Bach25DNEI,Bach21ONSA, Bach22WPFN,Barb76NSDE,Zeid90NFA2b,Roub13NPDE} and the references therein.

\subsection{Organization of the paper}
The paper is organized as follows: in Section 2, we set the analytical framework and briefly introduce some notions and results from the theory of convex analysis. In Section 3, we present and discuss the assumptions on the dissipation potential $\Psi$, the energy functional $\calE$, the perturbation $B$ as well as the external force $f$. Furthermore, we state the main result. Section 4 is devoted to the proof and in Section 5, we apply the abstract theory developed here to physically relevant examples which includes a mathematical model for visco-elasto-plasticity. In the Appendix, we collect certain  results from the subdifferential calculus.

\subsection{Notation and preliminaries}
For a proper functional $f:X\rightarrow (-\infty,+\infty]$ on a \textsc{Banach} space $(X,\Vert \cdot \Vert_X)$, we denote with the multivalued map $\partial f:X\rightrightarrows X^*$, the $($\textsc{Fr\'{e}chet}$)$ subdifferential of $f$ defined by
\begin{align*}
\partial f(u):=\left\lbrace \xi\in X^*: \liminf_{v\rightarrow u} \frac{f(v)-f(u)-\langle \xi,v-u\rangle_{X^*\times X}}{\Vert v-u\Vert_X} \geq 0\right \rbrace,
\end{align*} where $\langle\cdot,\cdot\rangle_{X^*\times X}$ denotes the duality pairing between the \textsc{Banach} space $X$ and its topological dual space $X^*$ equipped with the dual norm $\Vert \cdot \Vert_{X^*}:=\sup_{v\in X\backslash \lbrace 0\rbrace }\frac{\langle\cdot,v\rangle_{X^*\times X}}{\Vert v\Vert_X}$. The elements of the subdifferential are also called subgradients. If the set of subgradients of $f$ at a given point $u$ is nonempty, we say that $f$ is subdifferentiable at $u$. The effective domain of $f$ and the domain of its subdifferential $\partial f$ are defined and denoted by $\DOM(f):=\lbrace v\in X \mid f(v)<+\infty \rbrace$ and $\DOM(\partial f):=\lbrace v\in X : \partial f(v)\neq \emptyset\rbrace$, respectively.\\\\
Furthermore, we recall an important tool from the theory of convex analysis. For a
proper, lower semicontinuous, and convex functional $f: X\rightarrow (-\infty,+\infty]$, we define the so-called \textsc{Legendre--Fenchel} transform (or convex conjugate)  $f^*:X^*\rightarrow (-\infty,+\infty]$ by 
\begin{align*}
  f^*(\xi):=\sup_{u\in X}\left \lbrace \langle \xi,u\rangle_{X^*\times X} -
    f(u)\right \rbrace, \quad \xi\in X^*.
\end{align*} By definition, we directly obtain the \textsc{Fenchel--Young}
inequality
\begin{align*}
 \langle \xi,u\rangle_{X^*\times X} \leq  f(u)+f^*(\xi), \quad \text{for all } u\in X, \xi\in X^*.
\end{align*} It is easily checked that the transform itself is proper, lower semicontinuous and convex, see, e.g., \cite[Section 4, pp. 16]{EkeTem76CAVP}. If,  in addition, we assume $f(0)=0$, then $f^*(0)=0$ holds as well.\\\\
We recall also the following fact: let $(X,\Vert \cdot \Vert_{X})$ and $(Y, \Vert \cdot \Vert_Y)$ \textsc{Banach} spaces such that both $X$ and $Y$ are continuously embedded into another \textsc{Banach} space $Z$, and such that $X\cap Y$, equipped with the norm $\Vert \cdot \Vert_{X\cap Y}= \Vert \cdot \Vert_{X}+\Vert \cdot \Vert_Y$, is dense in both $X$ and $Y$. Then, the space $X\cap Y$ becomes a \textsc{Banach} space itself. If, furthermore, $X$ and $Y$ are separable and reflexive Banach spaces, the dual space can be identified by $X^*+Y^*$ with the dual norm $\Vert \xi \Vert_{X^*+Y^*}=\inf_{\overset{\xi_1\in X^*, \xi_2\in Y^*}{\xi=\xi_1+\xi_2}}\max \lbrace \Vert \xi_1\Vert_{X^*}, \Vert\xi_2\Vert_{Y^*}\rbrace$, and the duality pairing between $X\cap Y$ and $X^*+Y^*$ is given by
\begin{align*}
\langle f,v\rangle_{(X^*+ Y^*)\times (X\cap Y)}=\langle f_1,v\rangle_{X^*\times X}+\langle f_2,v\rangle_{Y^*\times Y} 
\end{align*} for all $v\in X\cap Y$ and any decomposition  $f=f_1+f_2$ with $f_1\in X$ and $f_2\in Y$, see, e.g., \cite[Kapitel 1, \S  5]{GaGrZa74NOOD}.\\\\  Furthermore, it is easily shown that $\rmL^p(0,T;X)\cap \rmL^p(0,T;Y)=\rmL^p(0,T;X\cap Y)$ for any $p\in [1,+\infty]$, where the measurability follows from the \textsc{Pettis} theorem, see, e.g., \cite[Theorem 2, p. 42]{DieUhl77VEME}.  If $X$ is separable and reflexive, the spaces $\rmL^p(0,T;X)$ are also separable and reflexive for all $1<p<\infty$ and $\rmL^\infty(0,T;X)$ is the dual of the separable space $\rmL^1(0,T;X^*)$. Finally, if the continuous embedding $X\hookrightarrow Y$ holds, then
\begin{align*}
\langle f,v\rangle_{X^*\times X}=\langle f,v\rangle_{Y^*\times Y} \quad \text{whenever } v\in X \text{ and } f\in Y^*.
\end{align*} see, e.g, \cite[Remark 3, pp. 136]{Brez11FASS} and \cite[Kapitel 1, \S  5]{GaGrZa74NOOD}. 
\section{Topological assumptions and the main result}
\label{se:AssumpExistRe.2}
\subsection{Function space setting}
We assume that $(U,\Vert\cdot\Vert_U), (V,\Vert\cdot\Vert_V)$, $(W,\Vert\cdot\Vert_W)$ and $(\widetilde{W},\Vert\cdot\Vert_{\widetilde{W}})$ are real, reflexive, and separable \textsc{Banach} spaces such that $U\cap V$ is separable and reflexive. Furthermore, we assume that $(H,\vert\cdot\vert,(\cdot,\cdot))$ is a \textsc{Hilbert} space with norm $\vert
\cdot\vert$ induced by the inner product $(\cdot,\cdot)$. Then, we assume the following dense, continuous, and compact embeddings 
\begin{align*}
\begin{cases}
U\cap V \overset{d}{\hookrightarrow}  U\overset{c,d}{\hookrightarrow} \widetilde{W} \overset{d}{\hookrightarrow} H \cong H^*\overset{d}{\hookrightarrow} \widetilde{W}^*\overset{d}{\hookrightarrow} U^* \overset{d}{\hookrightarrow}V^*+U^* \\
 U\cap V \overset{d}{\hookrightarrow} V\overset{c,d}{\hookrightarrow} W \overset{d}{\hookrightarrow} H \cong H^*\overset{d}{\hookrightarrow} W^*\overset{d}{\hookrightarrow} V^*\overset{d}{\hookrightarrow}V^*+U^*,
 \end{cases}
\end{align*} where $c$ and $d$ mean that the embedding is compact and dense, respectively. Moreover, if the perturbation does not explicitly depend on $u$ or $u'$, then we do not assume $U\overset{c}{\hookrightarrow} \widetilde{W}$ or $V\overset{c}{\hookrightarrow} W$, respectively. We further assume $V \hookrightarrow \widetilde{W}$ if $\calE_t^2\neq 0$, see Condition \ref{eq:cond.E2.1}. We note that we neither assume $U\hookrightarrow V$ nor $V \hookrightarrow U$ implying that the functionals $\calE_t$ and $\Psi_u$ live on different spaces. Since in this case, the subdifferential of $\Psi_u$ is nonlinear, we refer to the inclusion \eqref{eq:I.2} in the given framework as \textit{nonlinearly damped inertial system} $(U,V,W,\WW,H,\calE,\Psi,B,f)$.
\subsection{Assumptions on the functionals and operators}
In this section, we collect all the assumptions for the energy functional $\calE_t$, the dissipation potential $\Psi_u$, the perturbation $B$, and the external force $f$, and discuss them subsequently. We start with the assumptions for the dissipation potential $\Psi_u$.
\begin{enumerate}[label=\textnormal{(\thesection.$\Uppsi$\alph*)}, leftmargin=3.2em] 
\item \label{eq:Psi2.1} \textbf{Dissipation potential.} For every $u\in \DOM(\calE_t)$, let $\Psi_u: V\rightarrow [0,+\infty)$ be a lower semicontinuous and convex functional with $\Psi_u(0)=0$ such that the mapping $(u,v)\mapsto \Psi_u(v)$ is $\mathscr{B}(U)\otimes \mathscr{B}(V)$-measurable.

\item \label{eq:Psi2.2} \textbf{Superlinearity.} The functional $\Psi_u$ satisfies the following growth condition: there exists a positive real number $q>1$ such that for all $R>0$ there exist positive constants $c_R,C_R>0$ such that for all $u\in U$ with $\sup_{t\in [0,T]}\calE_t(u)\leq R$, there holds
\begin{align}\label{eq:Psi2.growth}
c_R(\Vert v\Vert^q_{V}-1)\leq \Psi_u(v)\leq C_R(\Vert v\Vert^q_{V}+1) \quad \text{for all } v\in V.
\end{align} 

\item \label{eq:Psi2.4}\textbf{Lower semicontinuity of $\Psi_u+\Psi_u^*$.} Let $v_n\rightharpoonup v$ in $\rmL^{q}(0,T; V)$, $\eta_n\rightharpoonup \eta$ in $\rmL^{q^*}(0,T; V^*)$, and $u_n(t) \rightharpoonup u(t)$ in $U$ for all $t\in[0,T]$ as $n\rightarrow \infty$ with $\sup_{t\in [0,T],n\in \mathbb{N}}\calE_t(u_n(t))<+\infty$  such that $\eta_n(t)\in \partial \Psi_{u_n(t)}(v_n(t))$ a.e. in $(0,T)$ for all $n\in \mathbb{N}$. Then, there holds 
\begin{align*}
\int_0^T \left(\Psi_{u(t)}(v(t))+ \Psi^*_{u(t)}(\eta(t))\right)\dd t\leq \liminf_{n\rightarrow \infty}\int_0^T \left(\Psi_{u_n(t)}(v_n(t))+ \Psi^*_{u_n(t)}(\eta_n(t))\right)\dd t.
\end{align*}
For the solvability of problem \eqref{eq:I.2}, only the previous assumptions are required. If we additionally assume the uniform monotonicity of $\partial\Psi_u$, we obtain stronger convergence of the discrete time-derivatives $\Von$ in the space $\rmL^q(0,T;V)$, see Lemma \ref{le:LimitPass2}.

\item \label{eq:Psi2.3} \textbf{Uniform monotonicity of $\partial\Psi_u$.} For all $R>0$, there exists a constant  $\mu_R>0$ such that 
\begin{align*}
\langle \xi-\eta,v-w\rangle_{V^*\times V} \geq \mu_R \Vert v-w\Vert_V^{\max\lbrace 2,q\rbrace}
\end{align*} for all $\xi\in \partial \Psi_u(v),\eta\in \partial  \Psi_u(w) $ and $u,v,w\in J_R:= \lbrace\tilde{v}\in V : \sup_{t\in [0,T]}\calE_t(\tilde{v})\leq R\rbrace$, where $q>1$ is from \eqref{eq:Psi2.growth}.
\end{enumerate}
\begin{rem}\label{re:Assump.Psi2} \mbox{} \vspace{-0.6em}
\begin{itemize}
\item[$i)$] We recall that the conjugate $\Psi_u^*:V^*\rightarrow \mathbb{R}$ is lower semicontinuous and convex itself, and that the growth condition \eqref{eq:Psi2.growth} implies the following growth condition for the conjugate $ \Psi^*_u$: for all $R>0$, there exist positive numbers $\bar{c}_R,\bar{C}_R>0$ such that for all $u\in U$ with $\sup_{t\in [0,T]}\calE_t(u)\leq R$, there holds
\begin{align*}
\bar{c}_R(\Vert \xi\Vert^{q^*}_{V^*}-1)\leq \Psi_u^*(\xi)\leq \bar{C}_R(\Vert \xi\Vert^{q^*}_{V^*}+1) \quad \text{for all }\xi\in V^*,
\end{align*} where $q^*=q/(q-1)$.
\item[$ii)$] It has been shown in \textsc{Stefanelli} \cite[Lemma 4.1]{Stef08BEPD} that the following stronger convergence in the sense of \textsc{Mosco} (we write $\Psi_{u_n}\Mto \Psi_u$) implies Condition \ref{eq:Psi2.4}: Let $u_n\rightharpoonup u\in V$ as $n\rightarrow \infty$. Then, for all $v\in V$, there holds
\begin{align} \label{Mosco}
\begin{cases}
a)\quad  \Psi_{u}(v)\leq \liminf_{n\to \infty} \Psi_{u_n}(v_n) \quad \text{for all }v_n\rightharpoonup v \text{ in } V,\\
b)\quad \exists \hat{v}_n\rightarrow v \text{ in $V$ such that } \Psi_u(v) \geq \limsup_{n\to \infty} \Psi_{u_n}(\hat{v}_n).
\end{cases}
\end{align}
\end{itemize}
\end{rem} 

Now, we proceed with the assumptions for the energy functional.
 
 \begin{enumerate}[label=\textnormal{(\thesection.E\alph*)},
 leftmargin=3.2em]  \label{eq:cond.E2}
\item \label{eq:cond.E2.1} 
\textbf{Basic properties.} For all $t\in[0,T]$, the
  functional $\calE_t:U \rightarrow \mathbb{R}$ is the sum of functionals $\calE^1:U\rightarrow \mathbb{R}$ and $\calE_t^2:\WW \rightarrow \mathbb{R}$. The functional $\calE^1(\cdot)=\frac{1}{2} b(\cdot,\cdot)$ is induced by a bounded, symmetric, and strongly positive bilinear form $b: U \times U\rightarrow \mathbb{R}$, i.e., there exist constants  $\mu,\alpha> 0$ such that
  \begin{align*}
   b(u,v)&\leq \alpha \Vert u\Vert_U\Vert v\Vert_U \quad \text{for all } u,v\in U\\
   \mu\Vert u\Vert_U^2 &\leq b(u,u) \qquad \quad \, \text{for all } u\in U.
    \end{align*} 
\item \label{eq:cond.E2.2}
 \textbf{Bounded from below.} $\calE_t$ is bounded from below uniformly in time, i.e., there exists a constant $C_0\in \mathbb{R}$ such that 
  \begin{align*}
  \calE_t(u)\geq C_0 \quad \text{for all } u\in U\text{ and }t\in[0,T].
  \end{align*} Since a potential is uniquely determined up to a constant, we assume without loss of generality $C_0=0$.
\item \label{eq:cond.E2.3} \textbf{Coercivity.} For every $t\in[0,T]$, $\calE_t$ has bounded sublevel sets in $U$, i.e., the set $J_{\lambda}=\lbrace u\in U : \sup_{t\in [0,T]}\calE_t(u)\leq \lambda \rbrace$ is bounded for all $\lambda\in \mathbb{R}$.
\item \label{eq:cond.E2.4} \textbf{Control of the time derivative.} For all $u\in \WW$,
  the mapping $t\mapsto \calE^2_t(u)$ is in $\rmC([0,T])\cap \rmC^1(0,T)$ and its derivative $\partial_t\calE^2_t$ is
  controlled by the function $\calE^2_t$, i.e., there exists $C_1>0$ such that
\begin{align}\label{eq:ass.E2.4}
  \vert \partial_t \calE^2_t(u)\vert \leq C_1 \calE^2_t(u)\quad \text{for all } t\in
  (0,T) \text{ and } u\in \WW^*.
\end{align} 
Furthermore, for all sequences $(u_n)_{n\in \mathbb{N}}\subset U$ with $u_n\rightharpoonup u$ in $U$ as $n\rightarrow \infty$ and $\sup_{n\in \mathtt{N}, t \in [0,T]}\calE_t(u_n)<+\infty$, there holds
\begin{align*}
\limsup_{n\rightarrow \infty}\partial_t\calE^2_{t}(u_n)\leq \partial_t\calE_t^2(u) \quad \text{for a.e. }t\in (0,T).
\end{align*}
\item \label{eq:cond.E2.6} \textbf{\textsc{Fr\'{e}chet} differentiability.} For all $t\in[0,T]$, the mapping $u\mapsto \calE^2_t(u)$ is \textsc{Fr\'{e}chet} differentiable on $\WW$ with derivative $\rmD\calE_t^2$ such that the mapping $(t,u)\mapsto \rmD\calE^2_t(u)$ is continuous as a mapping from $[0,T]\times \WW$ to $U^*$ on sublevel sets of the energy functional, i.e., for all $R>0$ and sequences $(u_n)_{n\in \mathbb{N}} \subset \WW $ and $(t_n)_{n\in \mathbb{N}} \subset [0,T]$ with $\sup_{t\in [0,T],n\in \mathbb{N}}\calE_t(u_n)<+\infty$,  $u_n\rightarrow u$ in $\WW$, and $t_n\rightarrow t\in [0,T]$ as $n\rightarrow \infty$, there holds 
\begin{align*}
\lim_{n\rightarrow \infty}\Vert \rmD\calE_{t_n}^2(u_n)-\rmD\calE_t^2(u)\Vert_{U^*}=0.
\end{align*} 
\item \label{eq:cond.E2.7}
 \textbf{$\lambda$-convexity.} 
There exists a non-negative real number $\lambda\geq 0$ such that
\begin{align*}
\calE_t(\vartheta u+(1-\vartheta )v)\leq &\vartheta \calE_t(u)+(1- \vartheta)\calE_t(v)+\vartheta(1-\vartheta)\lambda\vert u-v\vert^2
\end{align*} for all $t\in[0,T],  \vartheta\in [0,1]$, and $u,v\in U$.
\item \label{eq:cond.E2.8} \textbf{Control of $\rmD \calE^2_t$.}
  There exist positive constants $C_2>0$ and $\sigma>0$ such that
\begin{align*}
\Vert \rmD\calE^2_t(u) \Vert_{\WW^*}^\sigma \leq C_2 (1+\calE^2_t(u)+\Vert u\Vert_{\WW}) \quad \text{for all } t\in [0,T], u \in \WW.
\end{align*}
\end{enumerate}

Again, several remarks are in order. 

\begin{rem}\label{re:Assump.E2} \mbox{}\vspace{-0.6em}

\begin{itemize}
\item[$i)$] The assumptions on the quadratic form $\calE^1$ forces the \textsc{Banach} space $U$ to be a \textsc{Hilbert} space with the inner product induced by $b$. However, as we will see in Theorem \ref{th:MainExist2}, we can omit this assumption if the embedding $V \overset{c}{\hookrightarrow} \WW$ is compact, see Remark \ref{rem.b}.  Furthermore, the condition on the bilinear form $b$ implies that the \textsc{Fr\'{e}chet} derivative $\rmD\calE^1$ is given by a linear, bounded, symmetric and strongly positive operator $E\in \calL(U,U^*)$ such that $\calE^1(u)=\frac{1}{2}\langle \rmE u,u\rangle_{U^*\times U}$ is strongly convex and therefore sequentially weakly lower semicontinuous. Furthermore, the corresponding \textsc{Nemitski\v{i}} operator is a linear and bounded map from $\rmL^2(0,T;U)$ to $\rmL^2(0,T;U^*)$ and hence weak-to-weak continuous from $\rmL^2(0,T;U)$ to $\rmL^2(0,T;U^*)$. 
\item[$ii)$] From inequality \eqref{eq:ass.E2.4} in Assumption \ref{eq:cond.E2.4}, it follows after integration that 
\begin{align*}
  \begin{split}
  \sup_{t\in [0,T]}\calE_t^2(u)&\leq \ee^{C_1T} \inf_{t\in[0,T]} \calE_t^2(u), \\
 \vert \calE^2_t(u)-\calE^2_s(u)\vert &\leq e^{C_1T} \sup_{r\in [0,T]}\calE_r^2(u)\vert s-t\vert \quad \text{for all } u\in \WW^*, s,t\in [0,T].
 \end{split}
\end{align*} 
\item[$iii)$] The derivative of the $\lambda$-convex energy functional is characterized by the inequality
\begin{align}\label{eq:charact.subd}
\calE_t(u)-\calE_t(v)\leq \langle  \rmD\calE_t(u),u-v\rangle_{U^*\times U}+\lambda \vert u-v\vert^2
\end{align} for all $t\in [0,T]$, $u,v\in U$. In fact, the $\lambda$-convexity of $\calE_t$ can be replaced by inequality \eqref{eq:charact.subd}, since we only make use of \eqref{eq:charact.subd} in order to obtain a priori estimates, see Lemma \ref{le:DUEE2}. 
\end{itemize} 
\end{rem}
We recall that the \textsc{Fr\'{e}chet} differentiability of $\calE_t$ implies the subdifferentiability of $\calE_t$ and the subdifferential is a singleton given by $\partial \calE_t(u)=\lbrace \rmD \calE_t(u)\rbrace$.\\\\

Finally, we collect the assumptions concerning the perturbation $B$ and the external force $f$.

\begin{enumerate}[label= \textnormal{(\thesection.B\alph*)}, leftmargin=3.2em]  
\item \label{eq:B2.1} \textbf{Continuity.} The mapping $B:[0,T]\times \WW \times W \rightarrow V^*$ is continuous on  sublevel sets of $\calE_t$, i.e., for every converging sequence
$(t_n,u_n,v_n)\to (t,u,v)$ in $[0,T]\ti \WW \ti W$ as $n\rightarrow \infty$ with $\sup_{t\in [0,T], n\in\mathbb{N}}\calE_t(u_n)<+\infty$, there holds $B(t_n,u_n,v_n)\to B(t,u,v)$ in $V^*$ as $n\rightarrow \infty$.
\item \label{eq:B2.2} \textbf{Control of the growth.} There exist positive constants $\beta>0$ and $c,\nu \in(0,1)$  such that
\begin{align*}
c\,\Psi_u^*\left(\frac{-B(t,u,v)}{c}\right)\leq \beta (1+\calE_t(u)+\vert v\vert ^2+ \Psi_{u-rv}(v)^{\nu})
\end{align*} for all $u\in U,v\in V, \, t\in[0,T],$ and all $r\in [0,1)$.
\end{enumerate}

\begin{enumerate}[label= \textnormal{(\thesection.f)}, leftmargin=3.2em]  

\item \label{eq:f2} \textbf{External force.} There holds $f\in \rmL^{2}(0,T;H)$.
\end{enumerate}

\begin{rem}\label{re:Assump.f2} If the growth condition \ref{eq:Psi2.2} for $\Psi_u$ holds uniformly in $u\in U$, then more general external forces $f\in \rmL^1(0,T;H)+\rmL^{q^*}(0,T;V^*)$ can be considered, where $q^*>1$ is the conjugate number of $q$ from Condition \ref{eq:Psi2.2}.
\end{rem}
\subsection{Discussion of the assumptions}\label{su:Diss.NDS}
Apart from the remarks made above, we want to discuss certain conditions more in detail and provide concrete examples that satisfy the abstract setting. \\\\
As the name suggests, we consider evolution equations of second order with nonlinear damping, i.e., equations where $\partial\Psi_{u(t)}$ is nonlinear and in general multi-valued. As already mentioned in the literature review (Section \ref{lit.review}), this has not been studied before.\\\\
Ad Condition (\thesection.$\Uppsi$). The Condition \ref{eq:Psi2.1}  allows us to consider nonsmooth dissipation potentials. Furthermore, the assumption $\Psi_u(0)=0$ is not restrictive as the potential is uniquely determined up to a constant. The growth condition \ref{eq:Psi2.2} here is crucial to employ an integration by parts formula for the second derivative $u''$ proven in  \cite{EmmTha11DNEE}, see Lemma \ref{le:LimitPass2} below. Furthermore, as we mentioned in Remark \ref{re:Assump.Psi2} ii), the liminf estimate in Condition \ref{eq:Psi2.4} is already implied by the \textsc{Mosco}-convergence $\Psi_{u_n}\Mto \Psi_{u}$ for all sequences $u_n\rightharpoonup u$. The \textsc{Mosco}-convergence is related to the graph convergence of its subdifferential and stronger than the $\Gamma$-convergence \cite{Brai02GCB}. A prototypical example for a dissipation potential that fulfill Condition \ref{eq:Psi2.1}-\ref{eq:Psi2.3} is given by 
\begin{align*}
&\Psi_{\uu}(\vv)=\int_\Omega \left(  g_1(\nabla \uu)\frac{1}{p}\vert \nabla \vv(\xx)\vert^p+g_2(\nabla \uu)\vert \nabla \vv(\xx)\vert\right)\dd \xx \quad \text{or}\\
&\Psi_{\uu}(\vv)=\int_\Omega  \left( g_1(\uu)\frac{1}{p}\vert \vv(\xx)\vert^p+g_2(\uu)\vert  \vv(\xx)\vert\right)\dd \xx
\end{align*} on $V=\rmW_0^{1,p}(\Omega)^m$ or $V=\rmL^p(\Omega)^m$ with $m\in \mathbb{N}$ and $p\in (1,+\infty)$, respectively, where $g_1,g_2:\mathbb{R}^m\rightarrow \mathbb{R}$ are continuous function satisfying further conditions depending on the concrete form of the energy functional. See Chapter \ref{se:App}, where we discuss more general dissipation potentials. This type of dissipation potentials occur in rate-independent systems  such as in plasticity \cite{MieRou15RIST}, in ferromagnetic hysteresis \cite{Visi00FHO, MiRoSa13NADN} or in Elasto-visco-plasticity \cite{RajRou03EDSA}, see Section \ref{ex.martensitic}.
\\\\ 
Ad Condition (\thesection.E). The crucial assumption we make for the energy functional $\calE_t=\calE^1+\calE^2_t$ is that the leading part $\calE^1$ is defined by a bounded, symmetric, and strongly positive bilinear form $b:U\times U\rightarrow \mathbb{R}$. As mentioned in Remark \ref{re:Assump.E2}, the Conditions \ref{eq:cond.E2.4} and \ref{eq:cond.E2.6} ensure that we are able to control the time derivative of the energy functional and that we are able to pass to the limit in the time discretization scheme in the energy-dissipation inequality \eqref{sol:EDI2}, see Section \ref{se:Proof.2}. The Condition \ref{eq:cond.E2.8} is required in order to obtain bounds for the subgradients of $\calE_t$, which in turn is necessary to obtain a priori estimates for $u''$. The problem is that the bounds in Lemma \ref{le:DUEE2} only gives a priori estimates for the sum of the subgradient of $\calE_t$ and $u''$ which necessitates an independent bound that is given by Condition \ref{eq:cond.E2.8}. Condition \ref{eq:cond.E2.8} could be replaced by the more general assumption that $\partial\calE_t$ is a bounded operator.\\ A prototypical example for the energy functional is given by
\begin{align*}
\calE_t(\uu)=\int_\Omega \left( \frac{1}{2} \vert \nabla \uu\vert^2+\rmW(\uu)\right)\dd \xx+\int_\Omega \nabla \uu:\rmE \nabla \uu \dd \xx-\langle \mathbold{f}(t),\uu\rangle_{U^*\times U}
\end{align*} on $U=\rmH_0^{1}(\Omega)^m$ with $m\in \mathbb{N}$, where $\rmW:\mathbb{R} \rightarrow \mathbb{R}$ is a $\lambda$-convex and continuously differentiable function, e.g., $W(\uu)=(1-\uu^2)^2$, $\rmE:\mathbb{R}^m\rightarrow \mathbb{R}^m$ a uniformly positive definite and symmetric matrix, and $\mathbold{f}\in \rmC^1([0,T];U^*)$. This type of example occurs very often in models for ferro-magnetism where the solution $\uu$ is the so-called magnetization, see \cite{MiRoSa13NADN,RoMiSa08MACD}.\\\\
Ad Condition (\thesection.B). The continuity condition \ref{eq:B2.1} implies  that $B$ is a continuous perturbation of $\partial \calE_t$ and $\partial \Psi_u$. In practice, the term $B$ contains all non-variational and non-monotone contributions of lower order in terms of growth as well as spatial derivatives. This is reflected by Condition \ref{eq:B2.2}, where $B$ satisfies a growth condition in terms of the dissipation potential and its convex conjugate as well as the energy functional and the kinetic energy. In fact, the growth condition shows that the higher the order of the growths of $\Psi_u$ and $\calE_t$ are, the more we can allow for the growth of the perturbation. Condition \ref{eq:B2.2} ensures that we are able to control the growth of the perturbation in order to derive appropriate bounds. Both conditions can be generalized in a framework so that instead of a point-wise continuity and a pointwise growth condition, a continuity on suitable \textsc{Bochner} spaces can be imposed as well as a growth condition on the level of time integrals, see \cite{Akag08DNEB}. Furthermore, it would be sufficient to define the perturbation on the domain of the subdifferential of $\calE_t$, see, e.g.,  \cite{Akag08DNEB,Otan82NPNP}, where this has been considered for evolution inclusions of first order. A simple example for the perturbation is given by
\begin{align*}
B(t,\uu,\vv)=\int_\Omega \left( a(t) \vert \uu(\xx)\vert^p+ b(t)\vert \vv(\xx)\vert^q\right)\dd \xx 
\end{align*} on $W=\rmL^q(\Omega)^m$ and $\widetilde{W}=\rmL^q(\Omega)^m$ for appropriate $p,q\geq 1$, where $a,b \in \rmC(0,T)$. Obviously, $B$ is neither variational nor monotone.
\subsection{Statement of the main result}
Having discussed all assumptions, we are in the position to state the main result which includes the notion of a solution to \eqref{eq:I.2}.
\begin{thm}[Main result]
  \label{th:MainExist2} Let the nonlinearly damped inertial system \\ $(U,V,W,\widetilde{W},H,\calE,\Psi,B,f)$ be given and fulfill Assumptions \ref{eq:cond.E2.1}-\ref{eq:cond.E2.8}, \ref{eq:Psi2.1}-\ref{eq:Psi2.4} as well as \ref{eq:B2.1}-\ref{eq:B2.2} and \ref{eq:f2}. Then, for every $u_0\in U$ and $v_0\in H$, there exists a solution to \eqref{eq:I.2}, i.e., there exist functions $ u\in\rmC_w([0,T];U)\cap  \rmW^{1,\infty}(0,T;H)\cap \rmW^{2,q^*}(0,T;U^*+V^*)$ with $u-u_0\in \rmW^{1,q}(0,T;V)$ and $\eta \in \rmL^{q^*}(0,T;V^*)$ satisfying the initial conditions $u(0)=u_0$ in $U$ and $u'(0)=v_0$ in $H$ such that 
  \begin{align}\label{sol:IC2}
  \begin{split}
u''(t)+\eta(t)+\rmD \calE_t(u(t))+B(t,u(t),u'(t))=f(t) \quad &\text{in }U^*+V^*,\\
\eta(t)\in \partial \Psi_{u(t)}(u'(t)) \quad &\text{in } V^*,
\end{split}
\end{align} for almost every $t \in(0,T)$. Furthermore, the energy-dissipation balance 
\begin{align} \label{sol:EDI2}
\begin{split}
&\frac{1}{2}\vert u'(t) \vert^2 +\calE_t(u(t)) + \int_0^t \left(\Psi_{u(t)}(u'(r))+\Psi_{u(t)}^*(S(r)-\rmD \calE_r(r)-u''(r) \right)\dd r  \\
&= \frac{1}{2}\vert v_0 \vert^2 +\calE_0(u_0) + \int_0^t  \partial_r \calE_r(u(r)) \dd r+\int_0^t \langle S(r),u'(r) \rangle_{V^*\times V} \dd r
\end{split}
\end{align} holds for almost every $t\in(0,T)$, where $S(r):=f(r)-B(r,u(r),u'(r)),\, r\in[0,T]$. If $V \hookrightarrow U$, then \eqref{sol:EDI2} holds for all $t\in [0,T]$.
\end{thm}

\section{Proof of the main result}
The proof of Theorem \ref{th:MainExist2} consists of the following main steps:
\begin{itemize}
\item[\textbf{1)}] We discretize the inclusion in time via a semi-implicit \textsc{Euler} scheme with time step $\tau>0$ and show the solvability of the discrete problem.
\item[\textbf{2)}] We define interpolations functions and show that they satisfy a priori estimates.
\item[\textbf{3)}] We show compactness of the interpolation functions in suitable function spaces. 
\item[\textbf{4)}] We pass to the limit with $\tau\searrow 0$ and show existence of solutions that satisfy \eqref{sol:IC2} and \eqref{sol:EDI2}.
\end{itemize} The main difficulty in the proof consists of identifying the weak limits associated with the terms $\rmD \calE_t$ and $\partial\Psi_u$ in step \textbf{4)} where we have to derive and employ several new techniques.\\\\ 
In the following, each step will be carried out in a  subsection.
\subsection{Variational approximation scheme}\label{se:VarApprox.2}
The proof of Theorem \ref{th:MainExist2} relies on a semi-implicit time discretization scheme. Thus, for $N\in\mathbb{N}\backslash \lbrace 0\rbrace$, let
\begin{align*}
 I_\tau=\lbrace 0=t_0<t_1<\cdots< t_n=n\tau<\cdots <t_N=T \rbrace
\end{align*} be an equidistant partition of the time interval $[0,T]$ with step size $\tau:=T/N$, where we omit for simplicity the dependence of the nodes from the partition on the step size. Discretizing inclusion \eqref{eq:I.2} in a semi-implicit manner yields
\begin{align}
\label{eq:EuLa2}
\frac{\va-\vb}{\tau}+\partial \Psi_{\ub} \left(\va\right)+\rmD \calE_{\ta}(\ua) +B\left(\ta,\ub,\vb\right) \ni f_\tau^n \quad \text{in }U^*+V^*
 \end{align} for $n=1,\dots,N$ with $\va=\frac{\ua-\ub}{\tau}$. Here, we make the key observation that \eqref{eq:EuLa2} is the \textsc{Euler--Lagrange} equation associated with the functional 
 \begin{align*}
     u\mapsto \Upphi(\tau,t_{n-1},\ub,\uc,B(\ta,\ub,\vb)-f_\tau^n;u)
 \end{align*} given by
 \begin{align*}
  \Upphi(r,t,v,w,\zeta;u)= \frac{1}{2r^2} \vert u-2v+w\vert^2+
  r\Psi_v \left(\frac{u-v}{r}\right)+\calE_{t+r}(u)-\langle \zeta,u\rangle_{V^*\times V}
\end{align*} for $ r\in \mathbb{R}^{>0},t\in[0,T)$ with $r+t\in[0,T]$,  $u\in U\cap V, v\in V, w\in H$
and $ \zeta\in V^*$. This key observation enables us to determine the value $\ua$ recursively by minimizing the functional $\Upphi$ by the direct methods from the direct method in the calculus of variations and to allow us nonsmooth functionals $\Psi_u$ and $\calE_t$. We end up with the following variational approximation scheme
 \begin{align}\label{eq:ApproxSc2}
\begin{cases}
\un \in U\cap V \text{ and }\vn \in V \text{ are given; whenever $\ue,\dots,\ub \in U\cap V$ are known,}\\
\text{find } \ua \in J_{\tau,\tb}(\ub,\uc;B(\ta,\ub,\vb)-f_\tau^n)
\end{cases}
\end{align} for $n=1,\dots, N$, where $U^{-1}_\tau=\un-\vn \tau$ and $J_{r,t}(v,w;\zeta):=\mathrm{argmin}_{u\in U\cap V}\Upphi(r,t,v,w,\zeta;u)$.  The solvability of the discrete problem \eqref{eq:ApproxSc2} and that every solution fulfills the \textsc{Euler--Lagrange} equation \eqref{eq:EuLa2} is ensured by the following lemma.

\begin{lem}\label{le:Exist.Min2} Let the nonlinearly damped inertial system $(U,V,W,\widetilde{W},H,\calE,\Psi,B,f)$ be given and fulfill the Conditions \textnormal{ \ref{eq:cond.E2.1}-\ref{eq:cond.E2.3}},\textnormal{ \ref{eq:cond.E2.6}},\textnormal{ \ref{eq:cond.E2.7}}, and \textnormal{\ref{eq:Psi2.1}}-\textnormal{\ref{eq:Psi2.2}}. Furthermore, let $r\in (0,T)$ and $t\in [0,T)$ with $r+t\leq T$ as well as $v\in V, w\in H$ and $\zeta \in V^*$. Then, for all $r\leq \frac{1}{2\lambda}$, the set $J_{r,t}(v,w;\eta)$ is non-empty and single valued, where $\lambda$ is from \textnormal{\ref{eq:cond.E2.7}}. Furthermore, to every $u\in J_{r,t}(v,w;\zeta)$, there exists $\eta\in \partial \Psi_v\left(\frac{u-v}{r}\right)\subset V^*$ such that
  
  \begin{align*}
 \frac{u-2v-w}{r^2}+\eta+\rmD\calE_t(u)+\zeta=0 \quad \text{ in  }U^*+V^*.
 \end{align*}
  
\end{lem}
\begin{proof}[Proof] The proof follows from the direct methods of the calculus of variations as well as Lemma \ref{le:Subdif}.
Let $u\in U\cap V,v,w\in V, \zeta \in V^*$, and $r\in (0,\tau_0), t\in [0,T)$ with $r+t\leq T$ be given.  First of all, the
  \textsc{Fenchel--Young} inequality and the boundedness of the energy from below yield
\begin{align}
\label{eq:II.11}
   \Upphi(r,t,v,w,\zeta;u)&= \frac{1}{2r^2} \vert u-2v+w\vert^2+ r\Psi_v\left(\frac{u-v}{r}\right)+\calE_{t+r}(u)-\langle \zeta,u\rangle_{V^*\times V}\notag \\
  &\geq \frac{1}{2r^2} \vert u-2v+w\vert^2 -r\Psi_v^{*}(\zeta) +\calE_{t+r}(u)-\langle \zeta,v\rangle_{V^*\times V} \notag \\
    &\geq   \frac{1}{2r^2} \vert u-2v+w\vert^2-r\Psi_v^{*}(\zeta) -\langle \zeta,v\rangle_{V^*\times V}
\end{align} 
which implies on the one hand $\inf_{u\in U\cap V}\Upphi(r,t,v,w,\zeta;u)>-\infty$. On the other hand, we observe that
\begin{align}
\label{eq:II.12}
 \inf_{u\in U\cap V}\Upphi(r,t,v,w,\zeta;u)\leq \frac{1}{2r^2} \vert u_0-2v+w\vert^2+
  r\Psi_{v}\left(\frac{
 u_0-v}{r}\right)+\calE_{t+r}(u_0)-\langle \zeta,u_0\rangle_{V^*\times V}
\end{align} for any $u_0 \in U\cap V$,
so that $ \inf_{u\in U\cap V}\Upphi(r,t,v,w,\zeta;u)<+\infty$ holds as well. It remains to show that the global minimum is achieved by an element in $U\cap V$. In order to show this, let $(u_n)_{n \in
  \mathbb{N}}\subset U\cap V$ be a minimizing sequence for
$\Upphi(r,t,v,w,\zeta;\cdot)$. From \eqref{eq:II.11}, we deduce that
$(u_n)_{n \in \mathbb{N}}\subset U\cap V$ is contained in a sublevel set of $\calE_{t+r}$ and thus by Assumptions \ref{eq:cond.E2.3} and \ref{eq:Psi2.2} it is bounded in $U\cap V$. Since $U\cap V$ is a reflexive and separable \textsc{Banach} space, there exists
a subsequence (not relabelled) which converges weakly in $U\cap V$ towards a
limit $\tilde{u}\in U\cap V$. By the weak lower semicontinuity of the mapping
$u\mapsto \Upphi(r,t,v,w,\zeta;u)$, we have
\begin{align*}
  \Upphi(r,t,v,w,\zeta;\tilde{u})\leq \liminf_{n\rightarrow
    \infty}\Upphi(r,t,v,w,\zeta;u_n)=\inf_{\tilde{v}\in V}
  \Upphi(r,t,v,w,\zeta;\tilde{v}),
\end{align*} and therefore $u\in J_{r,t}(v,w;\zeta)\neq \emptyset$.
\end{proof}
Thus, Lemma \ref{le:Subdif} ensures that the minimizer of the mapping
\begin{align*}
u\mapsto \Upphi(\tau,\tb,\ub,\uc,B(\ta,\ub,\vb)-f_\tau^n;u),
\end{align*} fulfill the \textsc{Euler--Lagrange} equation \eqref{eq:EuLa2} for some $\eta\in \partial_{U\cap V}\Psi_v\left(\frac{u-v}{r}\right)\subset U^*+V^*$ where the subdifferential is taken on the space $U\cap V$ which can be realized by restricting the functional $\Psi_v$ to the space $U\cap V$. It remains to show that $\eta\in V^*$. Applying Lemma \ref{le:Leg.Fen}, there holds
\begin{align}\label{min.LF}
r\Psi_v\left(\frac{u-v}{r}\right)+r\widetilde{\Psi}^*_v\left(\eta \right)=\left \langle \eta, u-v\right \rangle_{(U^*+V^*)\times (U\cap V)}
\end{align} where $\widetilde{\Psi}^*_v$ is the convex conjugate of $\Psi_v$ on $U\cap V$. Taking into account $f=f^{**}$ for any proper, convex and lower semicontinuous functional (see \cite[Proposition 4.1, p. 18]{EkeTem76CAVP}) and equality \eqref{min.LF}, it is easy to see that
\begin{align*}
\widetilde{\Psi}^*_v(\xi)=
\begin{cases}
\Psi^*_v(\xi) \quad \text{if }\xi \in V^*\\
+\infty \quad \text{ else}
\end{cases}, \quad \xi \in U^*+V^*,
\end{align*} which immediately shows $\eta\in V^*$.
\subsection{Discrete Energy-Dissipation inequality and a priori estimates}\label{se:TimeDiscret.2}
In this section, we derive a priori estimates of the approximate solutions. Thus, let the initial values $u_0\in U\cap V$ and $v_0 \in V$ as well as the time step $\tau>0$ be given and fixed. As seen in Theorem \ref{th:MainExist2}, we will assume more general initial values $u_0\in U$  and $v_0\in H$ in the main existence result. These initial values will in the main proof be approximated by suitable sequences of values from $U\cap V$ and $V$, respectively. For now, let $\un:=u_0$ and $V_\tau^0=v_0$. Then, by Lemma \ref{le:Exist.Min2}, the variational approximation scheme \eqref{eq:ApproxSc2} is solvable and we obtain approximate values $(\ua)_{n=0}^N$ with $\un:=u_0$ and $V_\tau^0=v_0$. Next, from these approximate values, we define the piecewise constant and linear interpolations. The piecewise constant and linear interpolations are defined by
\begin{align}
\label{eq:Approx.U}
  &\overline{U}_\tau(0)=\underline{U}_\tau(0)=\widehat{U}_\tau(0):=U_\tau^0=u_0\, \text{ and } \notag \\
  &\underline{U}_\tau(t):=U_\tau^{n-1}, \quad \widehat{U}_\tau(t): =\frac{t_n-t}{\tau}U_\tau^{n-1}+\frac{t-t_{n-1}}{\tau}U_\tau^{n} \quad \text{for } t\in[t_{n-1},t_n),\\
  &\overline{U}_\tau(t):=U_\tau^n \quad \text{ for } t\in(t_{n-1},t_n]
  \quad \text{and}\quad \underline{U}_\tau(T)=U_\tau^N,\, n=1,\dots,N,\notag
\end{align} 
as well as
\begin{align}
\label{eq:Approx.V}
  &\overline{V}_\tau(0)=\underline{V}_\tau(0)=\widehat{V}_\tau(0):=V_\tau^0=v_0 \,\text{ and } \notag \\
  &\underline{V}_\tau(t):=V_\tau^{n-1}, \quad \widehat{V}_\tau(t): =\frac{t_n-t}{\tau}V_\tau^{n-1}+\frac{t-t_{n-1}}{\tau}V_\tau^{n} \quad \text{for } t\in[t_{n-1},t_n),\\
  &\overline{V}_\tau(t):=V_\tau^n \quad \text{ for } t\in(t_{n-1},t_n]
  \quad \text{and} \quad \underline{V}_\tau(T)=V_\tau^N,\, \ n=1,\dots,N, \notag
\end{align}  where $\va=\frac{\ua-\ub}{\tau}$ for $n=1,\dots,N$. We note that $\widehat{U}_\tau'=\overline{V}_\tau$ in the weak sense.
Furthermore, we define the piecewise constant function $f_\tau:[0,T]\rightarrow H$ by
\begin{align} \label{eq:Approx.xif}
&f_\tau (t)=f_\tau^n=\frac{1}{\tau}\int_{t_{n-1}}^{t_n}f(\sigma)\dd \sigma \quad \quad \text{ for } t\in[t_{n-1},t_n),\, n=1,\dots,N, \\
&f_\tau (T)=f_\tau^N.\notag
\end{align} 
Furthermore, by Lemma \ref{le:Exist.Min2}, there exists a sequence $(\eta_\tau^n)_{n=1}^N \subset V^*$ of subgradients fulfilling $\eta_\tau^n \in \partial \Psi_{U_\tau^{n-1}}(V_\tau^{n}),\, n=1,\dots, N$, such that
\begin{align*}
\frac{\va-\vb}{\tau}+\eta_\tau^n+\rmD \calE_{\ta}(\ua) +B\left(\ta,\ub,\vb\right) = f_\tau^n \quad \text{in }U^*+V^*, n=1,\dots,N.
\end{align*} Then, we define the measurable function $\eta_\tau:[0,T]\rightarrow V^*$ by 
\begin{align} \label{eq:Approx.xif2}
\eta_\tau(t)= \eta^n_{\tau} \quad \text{ for } t\in[t_{n-1},t_n),\, n=1,\dots,N, \quad \text{and}\quad \eta_\tau(T)= \eta^N_{\tau}.
\end{align} For notational convenience, we also introduce the piecewise constant functions $\teu_\tau:[0,T]\rightarrow [0,T]$ and $\teo_\tau:[0,T]\rightarrow [0,T]$ given by
\begin{align}\label{eq:Approx.t}
\begin{split}
&\teo_{{\tau}}(0): = 0\,\,\,\, \text{ and } \, \teo_{{\tau}}(t): =
t_n \quad \text{ for } t\in (t_{n-1},t_n], \\
&\teu_{{\tau}}(T): = T \, \text{ and } \, \teu_{{\tau}}(t): =
t_{n} \quad \text{ for } t\in [t_{n-1},t_n), \quad n=1,\dots,N.
\end{split}
\end{align} Obviously, there holds $\teo_\tau(t)\rightarrow t$ and
$\teu_\tau(t)\rightarrow t$ as $\tau\rightarrow 0$. \\\\ Furthermore, we introduce the short-hand notation 
\begin{align*}
    B_\tau(r):= B(\teo_\tau(r),\Uu (r),\Vu(r)) \quad \text{and} \quad S_\tau(r):= f_\tau(r)-B_\tau(r),\quad r\in[0,T].
\end{align*}
Having defined the interpolations, we are in the position to show a priori estimates in the following lemma.

\begin{lem}[A priori estimates]\label{le:DUEE2} Let the nonlinearly damped inertial system \\ $(U,V,W,\widetilde{W},H,\calE,\Psi,B,f)$ be given and fulfill Assumptions \ref{eq:cond.E2.1}-\ref{eq:cond.E2.8}, \ref{eq:Psi2.1}-\ref{eq:Psi2.4} as well as \ref{eq:B2.1}-\ref{eq:B2.2} and \ref{eq:f2}.  Furthermore, let $\overline{U}_\tau,\underline{U}_\tau, \widehat{U}_\tau,\overline{V}_\tau,\underline{V}_\tau, \widehat{V}_\tau, \eta_\tau$ and $f_\tau$ be the interpolations associated with the given values $u_0\in U\cap V$ and $v_0\in V$ as well as the step size $\tau>0$. Then, the discrete energy-dissipation inequality
\begin{align}
\label{eq: DUEE2}
 & \int_{\teo_\tau(s)}^{\teo_\tau(t)}\left(
    \Psi_{\Uu (r)} ( \Vo (r) ) +
    \Psi_{\Uu (r)}^*\left(S_\tau(r)-\Vh'(r)-
      \rmD\calE_{\teo(r)}(\Uo(r))\right) \right) \dd r
   \notag \\
   &\quad+\frac{1}{2} \left \vert \Vo(t)\right \vert^2+ \calE_{\teo_\tau(t)}(\overline{U}_\tau(t)) \notag \\ 
  &\leq \frac{1}{2} \left \vert \Vo(s) \right \vert^2+
  \calE_{\teo_\tau(s)}(\overline{U}_\tau(s))+\int_{\teo_\tau(s)}^{\teo_\tau(t)} \partial_r
  \calE_r(\Uu (r))\dd r+\int_{\teo_\tau(s)}^{\teo_\tau(t)}
  \langle S_\tau(r), \Vo(r) \rangle_{U^*\times U} \dd r\notag \\
 &\quad +\tau \lambda\int_{\teo_\tau(s)}^{\teo_\tau(t)}
     \vert \Vo(r)\vert^2\dd r
\end{align} holds for all $0\leq s< t\leq T$. Moreover, there exist positive constants $M,\tau^*>0$ such that the estimates 
\begin{align}
\label{eq:II2.29}
\sup_{t\in [0,T]}  \left \vert \Vo(t)\right \vert \leq M, \quad \sup_{t\in [0,T]} \calE_t(\Uo(t)) \leq M, \quad \sup_{t\in [0,T]}\vert \partial_t \calE_t(\Uu(t))\vert \leq M, \\ 
\label{eq:II2.30}
\int_0^T \left( \Psi_{\Uu (r)} \left( \Vo(r)\right) + \Psi_{\Uu (r)} ^*\left( S_\tau(r)-\Vh'(r)-  \rmD\calE_{\teo(r)}(\Uo(r))\right) \right) \dd r\leq M
\end{align} hold for all $0<\tau\leq \tau^*$. In particular, the families of functions 
\begin{subequations}
\label{eq:allbounds2}
\begin{align}
\label{eq:bUo2}
&(\Uo)_{0<\tau\leq \tau^*} \subset \rmL^\infty(0,T;U),\\
\label{eq:bVo2}
&(\Vo)_{0<\tau\leq \tau^*} \subset \rmL^q(0,T;V),\\
\label{eq:eta2}
&(\eta_\tau)_{0<\tau\leq \tau^*}\subset \rmL^{q^*}(0,T;V^*),\\
\label{eq:bVh'2}
&(\Vh')_{0<\tau\leq \tau^*}\subset \rmL^{\min\lbrace q^*,2\rbrace}(0,T;U^*+V^*),\\
\label{eq:bB2}
&(B_\tau)_{0<\tau\leq \tau^*}\subset \rmL^{\frac{q^*}{\nu}}(0,T;V^*),\\
\label{eq:bxi2}
&( \rmD\calE^2_{\teo}(\Uo))_{0<\tau\leq \tau^*} \subset \rmL^\infty(0,T;\WW^*),
\end{align} 
\end{subequations}
are uniformly bounded with respect to $\tau$ in the respective spaces, where $q^*>0$ is the conjugate exponent to $q>1$ and $\nu\in (0,1)$ being from Assumption \ref{eq:B2.2}. Finally, there holds
\begin{align}
\begin{split}
\label{eq:II2.31}
\sup_{t\in [0,T]}\left( \Vert \underline{U}_\tau(t)-\overline{U}_\tau(t)\Vert_V+\Vert \widehat{U}_\tau(t)-\overline{U}_\tau(t)\Vert_V\right)\rightarrow 0 \\
\sup_{t\in [0,T]}\left(\Vert \Vo(t)-\Vh(t)\Vert_{U^*+V^*} +\Vert \Vu(t)-\Vo(t)\Vert_{U^*+V^*}\right)
\rightarrow 0
\end{split}
\end{align} as $\tau \rightarrow 0$.
\end{lem}
\begin{proof}[Proof]
Let $(U_\tau^n)_{n=1}^N \subset U\cap V$ be the approximative values obtained from the variational approximation scheme \eqref{eq:ApproxSc2} which, by Lemma \ref{le:Subdif}, satisfy the \textsc{Euler--Lagrange} equation
\begin{align}\label{eq:dis.inc2}
 f_\tau^n -B(\ta,\ub,\vb) -\frac{\va-\vb}{\tau}-\rmD\calE_{t_n}(\ua) =\eta_\tau^n \in \partial \Psi_{\ub}(\va)
\end{align} for all $n=1,\dots,N$. According to Lemma \ref{le:Leg.Fen}, inclusion \eqref{eq:dis.inc2} is equivalent to
\begin{align}\label{estimate1}
&\Psi_{\ub} (\va)+\Psi_{\ub}^*\left(f_\tau^n -B(\ta,\ub,\vb) -\frac{\va-\vb}{\tau}-\rmD\calE_{t_n}(\ua)\right)\notag \\
&=\left \langle f_\tau^n -B(\ta,\ub,\vb) -\frac{\va-\vb}{\tau}-\rmD\calE_{t_n}(\ua),\va \right \rangle_{V^*\times V}, \quad n=1,\dots,N.
\end{align} Furthermore, the $\lambda$-convexity of $\calE_t$  \ref{eq:cond.E2.7} yields   
\begin{align}\label{estimate2}
-\left \langle \rmD\calE_{t_n}(\ua),\ua-\ub \right \rangle_{(U^*+V^*)\times (U\cap V)} &\leq \calE_{\ta}(\ub)-\calE_{\ta}(\ua)+\lambda \vert \ua-\ub \vert^2 \notag \\
&=\calE_{\tb}(\ub)-\calE_{\ta}(\ua)+\lambda \vert \ua-\ub \vert^2 \notag \\
&\quad +\int_{\tb}^{\ta} \partial_r \calE_r(\ub) \dd r
\end{align} for all $n=1,\dots,N$, see Remark \ref{re:Assump.E2} $iii)$.
Then, plugging in the inequality \eqref{estimate2} into \eqref{estimate1} and making use of the identity 
\begin{align} 
\left( \va-\vb,\va \right)= \frac{1}{2} \left( \vert \va \vert^2-\vert \vb \vert^2+ \vert \va-\vb \vert^2\right) \quad n=1,\dots,N,
\end{align} as well as the fact that $\langle w,v\rangle_{V^*\times V} =(w,v)$ whenever $v\in V$ and $w\in H$, we obtain
\begin{align*}
&\frac{1}{2} \vert \va \vert^2 +\calE_{\ta}(\ua)+ \tau \Psi_{\ub}(\va)+\tau\Psi_{\ub}^*\left(S_\tau^n -\frac{\va-\vb}{\tau}-\rmD\calE_{t_n}(\ua)\right)\notag \\
 &\leq \frac{1}{2} \vert \vb \vert^2 +\calE_{\tb}(\ub)+\int_{\tb}^{\ta} \partial_r \calE_r(\ub) \dd r +\lambda \int_{\tb}^{\ta} \vert \va \vert^2 \dd r+\tau\left \langle S_\tau^n,\va \right \rangle_{V^*\times V}
\end{align*} for all $n=1,\dots,N$, where $S_\tau^n:=f_\tau^n -B(\ta,\ub,\vb), n=1,\dots, N$. Summing up the inequalities over $n$ yields \eqref{eq: DUEE2}. The estimates \eqref{eq:II2.29} and \eqref{eq:II2.30} are obtained by employing the discrete version of \textsc{Gronwall}'s lemma and the following estimates. First, employing Condition \ref{eq:B2.2} and the \textsc{Fenchel--Young} inequality, we obtain for $\varepsilon<\frac{1-c}{2}$ 
\begin{align*}
\tau \langle S_\tau^n,\va \rangle&= \tau \langle -B(\ta,\ub,\vb)+f_\tau^n ,\va \rangle_{V^*\times V}\notag \\
&= \tau \langle -B(\ta,\ub,\vb),\va \rangle_{V^*\times V}+\tau \langle f_\tau^n,\va \rangle_{V^*\times V}\notag \\
&\leq  c \tau \Psi_{\ub}(\va)+c\tau\Psi_{\ub}^*\left(\frac{-B(\ta,\ub,\vb)}{c}\right)+ \frac{\tau}{2} (\vert f_\tau^n \vert^2+\vert \va \vert^2)\notag \\
&\leq c \tau \Psi_{\ub}(\va) + \tau \beta (1+\calE_{\ta}(\ub)+\vert \vb \vert^2 +\Psi_{\ub-r \vb}(\vb)^\nu)\\
&\quad+\frac{\tau}{2} (\vert f_\tau^n \vert^2+\vert \va \vert^2),\notag\\
&\leq c \tau \Psi_{\ub}(\va) + \tau \beta (1+\calE_{\ta}(\ub)+\vert \vb \vert^2)+\tau \varepsilon \Psi_{\ub-r \vb}(\vb)\\
&\quad +\frac{\tau}{2} (\vert f_\tau^n \vert^2+\vert \va \vert^2)+\tau C,\notag
\end{align*} for all $r\in [0,1]$ and a constant $C=C(\beta,\nu,\varepsilon)>0$. For $n=1$, we choose $r=0$ and for $n\geq 2$, we choose $r=\tau$ in Condition \ref{eq:B2.2} and note that $\ub-\tau \vb=U_\tau^{n-2}$. Second, using Condition \ref{eq:cond.E2.4}, there holds 
\begin{align*}
\int_{\tb}^{\ta} \partial_r \calE_r(\ub) \dd r\leq \int_{\tb}^{\ta} C_1 \calE_r(\ub) \dd r \leq C_1 \int_{\tb}^{\ta}  \sup_{t\in [0,T]}\calE_t(\ub) \dd r.
\end{align*} 
Inserting the obtained inequalities into \eqref{eq: DUEE2} and summing up all inequalities from $1$ to $n$, we find another positive constant $C>0$ such that
\begin{align*}
&\frac{1}{2} \vert \va \vert^2 +\frac{1}{C_1}\sup_{t\in [0,T]}\calE_t(\ua)+ \int_{0}^{\ta}\left((1- \alpha(\tau)) \Psi_{\Uu(r)}(\Vo(r))\right) \dd r \notag \\
&\quad + \int_{0}^{\ta} \Psi_{\Uu(r)}^*\left(S_\tau (r) -\Vh'(r)-\xi_{\tau}(r) \right) \dd r \notag \\
 &\leq C (1+\vert v_0 \vert^2 +\calE_{0}(u_0)+\Vert f \Vert^2_{\rmL^2(0,T;H)}+\Psi_{u_0}(v_0) + \int_{0}^{\ta} \big( \vert \Vo(r)\vert^2+ \sup_{t\in [0,T]}\calE_t(\Uo(r))\big)\dd r) 
\end{align*}  where $\alpha(\tau):=c+\ct+\tau\frac{\lambda}{\mu}<1$ for all $\tau<\tau^*: =\frac{\mu}{\lambda}(1-c-\tilde{c})$ and we made use of the estimate for the interpolation $f_\tau$
\begin{align}\label{eq:f.est}
\Vert f_\tau\Vert^2_{\rmL^2(0,T;H)}&=\sum_{k=1}^N \tau \vert f_\tau^k\vert^2\notag\\
&=\sum_{k=1}^N \frac{1}{\tau}\vert \int_{t_{k-1}}^{t_k}f(\sigma)\dd \sigma \vert^2\notag \\
&\leq  \sum_{k=1}^n \int_{t_{k-1}}^{t_k} \vert f(\sigma)\vert^2\dd \sigma= \int_0^{T}\vert f(\sigma)\vert^2\dd \sigma= \Vert f\Vert^2_{\rmL^2(0,T;H)}.
\end{align} Then, by the discrete version of \textsc{Gronwall}'s lemma (see, e.g., \cite[Lemma 3.2.4, p. 68]{AmGiSa05GFMS}), there exists a constant $M>0$ such that \eqref{eq:II2.29} and \eqref{eq:II2.30} are satisfied. Now, we seek to show the bounds in \eqref{eq:allbounds2}. Due to the bounds obtained in \eqref{eq:II2.29} and \eqref{eq:II2.30}, the coercivity of $\Psi_u$ and $\Psi_u^*$ yield the boundedness of $(\Vo)_{0<\tau\leq \tau^*}\subset \rmL^q(0,T;V)$ and $(\eta_\tau)_{0<\tau\leq \tau^*}=(S_\tau-\Vh'- \rmD\calE^2_{\teo}(\Uo)))_{0<\tau\leq \tau^*}\subset \rmL^{q^*}(0,T;V^*)$ uniformly in $\tau$. The uniform boundedness of $(B_\tau)_{0<\tau\leq \tau^*}\subset \rmL^{\frac{q^*}{\nu}}(0,T;V^*)$ in turn follows from Assumptions \ref{eq:B2.2} and \ref{eq:Psi2.2}: taking into account that for all $\zeta\in V^*$ the mapping $r\mapsto r\Psi^*(\zeta/r)$ is monotonically decreasing on $(0,+\infty)$ which follows from the convexity of $\Psi^*$ and $\Psi^*(0)=0$, we obtain another constant $C_B>0$ such that 
\begin{align}\label{eq:II.38.1}
\int_0^T \Vert B_\tau(r))\Vert_*^{\frac{q^*}{\nu}} \dd r&\leq \int_0^T C_B\left( \Psi_{\Uu(r)}^*\left( B(\teo_{\tau}(r),\Uu(r),\Vu(r))\right)^\frac{1}{\nu}+1\right)\dd r \notag\\
&\leq \int_0^TC_B\left( c\,\Psi_{\Uu(r)}^*\left(\frac{ B(\teo_{\tau}(r),\Uu(r),\Vu(r))}{c}\right)^\frac{1}{\nu}+1\right)\dd r \notag\\
&\leq \int_0^T C \left( 1+\calE_{\teu_{\tau}(r)}(\Uu(r))^\frac{1}{\nu}+\vert \Vu(r)\vert^\frac{2}{\nu}+ \Psi_{\Uu(r)}\left( \Vu(r) \right)\right)\dd r \notag\\
&\leq N 
\end{align} for positive constants $C,N>0$ independent of $\tau$, where $c\in (0,1)$ is from Assumption \ref{eq:B2.2}. Since $(f_\tau)_{0<\tau\leq \tau^*}$ is uniformly bounded in $\rmL^2(0,T;H)$, it follows  that $(\Vh'+\xi_\tau)_{0<\tau\leq \tau^*}$ is uniformly bounded in $\rmL^{\min\lbrace q^*,2\rbrace}(0,T;V^*)$. Finally, Assumption \ref{eq:cond.E2.8} implies a uniform bound for $(\rmD\calE^2_{\teo}(\Uo))_{0<\tau\leq \tau^*}$ in $\rmL^\infty(0,T;W^*)$. Since all previous families of functions are bounded in the common space $\rmL^{\min\lbrace q^*,2\rbrace}(0,T;U^*+V^*)$, we deduce that $(\Vh')_{0<\tau\leq \tau^*}$ is uniformly bounded in $\rmL^{q^*}(0,T;U^*+V^*)$ with respect to $\tau$.\\ Finally, the convergences \eqref{eq:II2.31} follow from the bounds of $(\Vh')_{0<\tau\leq \tau^*}\subset \rmL^{\min\lbrace q^*,2\rbrace}(0,T;U^*+V^*)$ and $(\Vo)_{0<\tau\leq \tau^*}\subset \rmL^q(0,T;V)$ as well as the estimates 
\begin{align*}
&\Vert \Uh(t)-\Uo(t)\Vert_V \leq \Vert \Uu(t)-\Uo(t)\Vert_V = \int_{\teu(t)}^{\teo(t)} \Vert \Vo(r)\Vert_V \dd r \quad \text{ and } \\
&\Vert \Vh(t)-\Vo(t)\Vert_{U^*+V^*} \leq \Vert \Vu(t)-\Vo(t)\Vert_{U^*+V^*}= \int_{\teu(t)}^{\teo(t)} \Vert \Vh'(r)\Vert_{U^*+V^*}\dd r 
\end{align*}
 for all $t\in[0,T]$ which completes the proof.
\end{proof} 
\subsection{Compactness}
\label{su:compactness2}
In this section, we prove the compactness of the interpolation functions  in suitable \textsc{Bochner} spaces with respect to a suitable topology. This will enable us to pass to the limit in the discrete energy-dissipation inequality \eqref{eq: DUEE2} and the discretized inclusion \eqref{eq:EuLa2} as $\tau\searrow 0$. After identifying all of the limits, we will indeed obtain a solution to the \textsc{Cauchy} problem \eqref{eq:I.2}. The compactness result is given in the following lemma.

\begin{lem}[Compactness] \label{le:LimitPass2} Under the same assumptions of Lemma
  \ref{le:DUEE2}, 
 let $(\tau_n)_{n \in \mathbb{N}}$ be a
  vanishing sequence of step sizes, and let $u_0\in U\cap V $ and $v_0\in V$. Then, there exists a subsequence, still denoted by $(\tau_n)_{n \in \mathbb{N}}$, a pair of functions $(u,\eta)$ with   \begin{align*}
  &u\in\rmC_{w}([0,T];U)\cap \rmW^{1,q}(0,T;V) \cap\rmW^{1,\infty}(0,T;H)\cap \rmW^{2,q^*}(0,T;U^*+V^*) \text{ and}\\
  &\eta\in \rmL^{q^*}(0,T;V^*),
\end{align*} fulfilling the initial values $u(0)=u_0$ in $U$ and $u'(0)=v_0$ in $H$ such that the following convergences hold
\begin{subequations}
\label{eq:LP2}
\begin{align}
\label{eq:LP.uhu.weak2}
\Uun,\Uon,\Uhn \overset{*}{\rightharpoonup} u \quad \text{in } &\rmL^{\infty}(0,T;U\cap V),\\
\label{eq:B.uhh2}
\Uon, \Uun \rightarrow u \quad \text{in } &\rmL^r(0,T;\WW) \quad \text{for any } r\geq 1,\\
\label{eq:LP.vo2}
\Von, \Vun \overset{*}{\rightharpoonup} u' \quad \text{in } &\rmL^q(0,T;V)\cap \rmL^{\infty}(0,T;H),\\
\label{eq:LP.vo.strong2}
\Von, \Vun\rightarrow u' \quad \text{in } &\rmL^q(0,T;W)\\
\label{eq:LP.eta}
\eta_{\tau_n} \rightharpoonup \eta \quad \text{in } &\rmL^{q^*}(0,T;V^*),\\
\label{eq:LP.vhn'2}
\Vhn' \rightharpoonup u'' \quad \text{in } &\rmL^{\min\lbrace 2,q^*\rbrace}(0,T;U^*+V^*),\\
\label{eq:LP.xi2}
\rmE \Uon \rightharpoonup \rmE u \quad \text{in } &\rmL^2(0,T;U^*),\\
\label{eq:LP.xi22}
\rmD\calE^2_{\teon}(\Uon) \rightarrow \rmD\calE^2_t(u) \quad \text{in } &\rmL^r(0,T;\WW^*)\quad \text{for any } r\geq 1,\\
\label{eq:LP.f2}
f_{\tau_n}\rightarrow f \quad \text{in } &\rmL^{2}(0,T;H),\\
\label{eq:LP.B2}
B_{\tau_n}\rightarrow B(\cdot,u(\cdot),u'(\cdot))\quad \text{in } &\rmL^{q^*}(0,T;V^*),\\
\label{eq:uhu.weak2}
\Uhn(t),\Uun(t),\Uon(t) \rightharpoonup u(t) \quad \text{in } &U \,\, \text{ for all }t\in[0,T],\\
\label{eq:uhu.weakinV2}
\Uun(t) \rightharpoonup u(t) \quad \text{in } &V \,\, \text{ for all }t\in[0,T],\\
\label{eq:B.uhh.ptw2}
\Uhn(t),\Uun(t),\Uon(t) \rightarrow u(t) \quad \text{in } &\WW \,\text{ for all }t\in[0,T],\\
\label{eq:vh.ptw.weak2}
\Vhn(t), \Vun(t), \Von(t) \rightharpoonup u'(t) \quad \text{in } & V \,\,\, \text{ for all }t\in[0,T],\\
\label{eq:vh.ptw.strong}
\Vhn(t), \Vun(t), \Von(t) \rightarrow u'(t) \quad \text{in } & W \,\, \text{ for all }t\in[0,T].
\end{align}
\end{subequations} 
Furthermore, if the dissipation potential also satisfies Assumption \ref{eq:Psi2.3}, then, there holds 
\begin{subequations}
\label{eq:LP2all}
\begin{align}
\label{eq:C.uhu.weak2.strong}
\Von \rightarrow u' \quad \text{in } &\rmL^{\max\lbrace 2,q\rbrace}(0,T;V),\\
\label{eq:LP.uhu.weak2.strong}
\Uhn \rightarrow u \quad \text{in } &\rmC([0,T];V).
\end{align} 
\end{subequations} Finally, the function $u$ satisfies the inequality
\begin{align}\label{formulaxy}
&\frac{1}{2}\vert v_0 \vert^2-\frac{1}{2}\vert u'(t) \vert^2+\calE_{0}(u_0)-\calE_{t}(u(t))+\int_{0}^{t}\partial_r\calE_r(u(r))\dd r \notag \\
&\leq -\int_0^{t}\langle u''(r)+\rmD\calE_r(u(r)),u'(r)\rangle_{V^*\times V}
\end{align} for almost every $t\in(0,T)$.
\end{lem}
\begin{proof}[Proof] 
The convergences \eqref{eq:LP.uhu.weak2}, \eqref{eq:LP.vo2}, \eqref{eq:LP.eta}, \eqref{eq:LP.vhn'2}, and \eqref{eq:LP.xi2} follow (up to a subsequence) from the bounds shown in \eqref{eq:bUo2}-\eqref{eq:eta2} and Remark \ref{re:Assump.E2} i). We note that by standard arguments, we can identify the weak limits in \eqref{eq:LP.vo2} and \eqref{eq:LP.vhn'2} with $u'$ and $u''$, respectively, denoting with $u$ the weak limit in \eqref{eq:LP.uhu.weak2}. We recall the fact that $\rmL^{\infty}(0,T;X)\cap \rmC_w([0,T];Y)=\rmC_w([0,T];X)$ for two \textsc{Banach} spaces $X$ and $Y$ with $X$ being reflexive such that the continuous and dense embedding $X\overset{d}{\hookrightarrow} Y$ holds, see, e.g., in \textsc{Lions \& Magenes} \cite[Lemma 8.1, p. 275]{LioMag72NHBV}. Based on this fact, we derive that for $X=U$ and $Y=H$, there holds $u\in \rmC_w([0,T];U)$.  The convergences \eqref{eq:B.uhh2} and \eqref{eq:LP.vo.strong2} follow from the \textsc{Lions--Aubin--Dubinski\v{i}} lemma\footnote{The \textsc{Lions--Aubin--Dubinski\v{i}} lemma is a discrete version of the classical \textsc{Lions--Aubin} lemma.} (see, e.g., \cite[Theorem 1]{DreJun12CFPC}). We proceed with proving the pointwise convergence \eqref{eq:uhu.weak2}-\eqref{eq:vh.ptw.strong}. First, we note that from $\Vhn\in \rmW^{1,1}(0,T;U^*+V^*)\hookrightarrow \rmC([0,T];U^*+V^*)$ and \eqref{eq:LP.vhn'2}, there holds $\Vhn(t)\rightharpoonup u'(t)$ in $U^*+V^*$ as $n\rightarrow \infty$ for all $t\in [0,T]$. Since $\Vhn(t)$ is uniformly bounded in $H$ for all $t\in [0,T]$, it is (up to a subsequence) weakly convergent to $u'(t)$ in $H$ for all $t\in [0,T]$. Since the weak limit is unique in $U^*+V^*$, we obtain with the subsequence principle the convergence of the whole sequence. Together with the strong convergence in \eqref{eq:II2.31}, this implies \eqref{eq:vh.ptw.weak2}. With the same argument, we can show  the pointwise weak convergences \eqref{eq:uhu.weak2} and \eqref{eq:uhu.weakinV2} where for \eqref{eq:uhu.weakinV2} we also use the fact that $u_0\in U\cap V$.\\ Since $\rmC([0,T];H)$ is dense in $\rmL^2(0,T;H)$, there exists for every $\epsilon>0$ a function $f^\varepsilon\in \rmC([0,T];H)$ such that $\Vert f^\varepsilon-f\Vert_{\rmL^2(0,T;H)}<\varepsilon/3$. We obtain for sufficiently small step sizes $\tau_n$
\begin{align*}
\Vert f_{\tau_n}-f\Vert_{\rmL^2(0,T;H)}&\leq \Vert f_{\tau_n}-f^\varepsilon_{\tau_n} \Vert_{\rmL^2(0,T;H)}+\Vert f^\varepsilon_{\tau_n}-f^\varepsilon\Vert_{\rmL^2(0,T;H)}+\Vert f^\varepsilon-f\Vert_{\rmL^2(0,T;H)}\\
&\leq\Vert f-f^\varepsilon \Vert_{\rmL^2(0,T;H)}+\Vert f^\varepsilon_{\tau_n}-f^\varepsilon\Vert_{\rmL^2(0,T;H)}+\Vert f^\varepsilon-f\Vert_{\rmL^2(0,T;H)}\\
&\leq \varepsilon/3+\varepsilon/3+\varepsilon/3=\varepsilon,
\end{align*} where we also used the estimate \eqref{eq:f.est} for the first term. The second term can be made smaller than $\varepsilon/3$ for sufficiently small step sizes because of the uniform continuity of $f^{\varepsilon}$. Now, we want to show the convergence of the perturbation in \eqref{eq:LP.B2}. To do so, we denote by $\calB(u)(t)=B(t,u(t),u'(t)),t\in[0,T]$, the associated \textsc{Nemitski\v{i}} operator and recall that $B_{\tau_n}(t)=B(\teon(t),\Uun(t),\Vun(t)), t\in[0,T]$. First, the pointwise convergences \eqref{eq:B.uhh.ptw2} and \eqref{eq:vh.ptw.strong} as well as the continuity condition \ref{eq:B2.1} imply
\begin{align}\label{eq:B.conv}
\Vert \calB_{\tau_n}(t)- \calB(u)(t) \Vert_{V^*} \rightarrow 0 \quad \text{as $n\rightarrow \infty$ a.e. in }(0,T).
\end{align} By the growth condition \ref{eq:B2.2}, we can show that $\calB(u)\in \rmL^\frac{q^*}{\nu}(0,T;V^*)$ in the same way as in \eqref{eq:II.38.1}. Hence, $B_{\tau_n}-\calB(u)\in \rmL^\frac{q^*}{\nu}(0,T;V^*)$ is uniformly bounded in $n\in \mathbb{N}$ by a constant $\widetilde{M}>0$. Then, by \textsc{Egorov}'s theorem, for every $\varepsilon>0$ there exists a subset $E\subset [0,T]$ with measure $\mu(E)<\varepsilon$ such that 
\begin{align*}
\lim_{n\rightarrow \infty}\sup_{t\in [0,T]\backslash E}\Vert B_{\tau_n}(t)- \calB(u)(t) \Vert_{V^*}=0.
\end{align*} Therefore, for every $\varepsilon>0$ there exists an index $\tilde{N}\in \mathbb{N}$ such that for all $n\geq \tilde{N}$, there holds
\begin{align*}
\Vert B_{\tau_n}(t)- \calB(u)(t)\Vert_{V^*}\leq \varepsilon \quad \text{for all } t\in [0,T]\backslash E.
\end{align*} Invoking the latter estimate, we obtain for all $n\geq \tilde{N}$
\begin{align*}
&\Vert B_{\tau_n}-\calB(u)\Vert_{\rmL^{q^*}(0,T;V^*)}\\
&\leq \left(  \int_E \Vert B_{\tau_n}(t)-\calB(u)(t) \Vert_{V^*}^{q^*}\dd t\right)^{\frac{1}{{q^*}}}+\left(\int_{[0,T]\backslash E} \Vert B_{\tau_n}(t)-\calB(u)(t) \Vert_{V^*}^{q^*}\dd t\right)^{\frac{1}{{q^*}}}\\
&\leq \mu(E)^{1-\nu} \left(  \int_E \Vert B_{\tau_n}(t)-\calB(u)(t)\Vert_{V^*}^\frac{q^*}{\nu}\dd t\right)^\frac{\nu}{q^*}+\varepsilon T^{\frac{1}{q^*}}\\
&\leq \varepsilon^{1-\nu}\tilde{M}+\varepsilon T^{\frac{1}{q^*}}
\end{align*} and hence \eqref{eq:LP.B2}. Further, from the growth condition \ref{eq:cond.E2.8}, we obtain for all $t\in [0,T]$
\begin{align*}
\Vert \rmD\calE_{\teon(t)}^2(\Uon(t))\Vert_{\WW^*}^\sigma\leq C_3(1+\calE^2(\Uon(t))+\Vert \Uon(t)\Vert_{\WW})
\end{align*} and in view of the a priori estimates \eqref{eq:II2.29},
\begin{align*}
\Vert \rmD\calE_{\teon(t)}^2(\Uon(t))\Vert_{\WW^*} &\leq C \quad \text{for all }t\in [0,T].
\end{align*} Together with the convergence  \eqref{eq:B.uhh.ptw2} and the continuity condition \ref{eq:cond.E2.6}, this leads to \eqref{eq:LP.xi22}. The  assertions \eqref{eq:C.uhu.weak2.strong} and \eqref{eq:LP.uhu.weak2.strong} follow immediately from Assumption \ref{eq:Psi2.3} and
\begin{align*}
&\limsup_{n\rightarrow \infty}\int_0^T \Vert \Von(r)-u'(r)\Vert^{\max\lbrace p,2\rbrace}_V \dd r\\
&\leq \limsup_{n\rightarrow \infty}\int_0^T \langle \eta_n(r)-\eta(r),\Von(r)-u'(r)\rangle_{V^*\times V}\dd r\leq 0
\end{align*} and $\eta(t)\in \partial \Psi_{u(t)}(u'(t))$ a.e. in $(0,T)$, which we will show in the proof of the main result. It remains to show the inequality \eqref{formulaxy}. The difficulty in proving the aforementioned inequality is that we are not allowed to split the duality pairing in the integral on the right-hand side and consider each integral separately as only the sum lives on the space $V^*$. However, we can circumvent this problem following the proof of \cite[Lemma 6]{EmmTha11DNEE}. The idea relies on the  regularization of the function $u'$ by its so-called \textsc{Steklov} average.  For a function $v\in \rmL^p(0,T;X), p\geq 1,$ defined on a \textsc{Banach} space $X$ and being extended by zero outside $[0,T]$, the \textsc{Steklov} average is, for sufficiently small $h>0$, given by 
\begin{align*}
S_h v(t):=\frac{1}{2h}\int_{t-h}^{t+h} v(r)\dd r.
\end{align*} It is readily seen that $S_hv\in \rmL^p(0,T;X)$ and $\Vert S_hv\Vert_{\rmL^p(0,T;X)}\leq \Vert v\Vert_{\rmL^p(0,T;X)}$. Furthermore, it can be shown by a regularization argument that $S_hv\rightarrow v$ in $\rmL^p(0,T;X)$ as $h\rightarrow 0$, see , e.g., \cite[Theorem 9, p. 49]{DieUhl77VEME}. \\ Defining $Kv(t)=\int_0^t v(r)\dd r$, we commence with calculating
\begin{align*}
 &\quad -\int_s^{t}\langle (S_h u')'(r)+\rmD\calE_r(u_0+(KS_h u')(r)),(S_h u')(r)\rangle_{V^*\times V}\dd r\\
 &=-\int_s^{t}\langle (S_h u')'(r)+\rmE(u_0+(KS_h u')(r))+\rmD\calE^2_r(u_0+(KS_h u')(r)),(S_h u')(r)\rangle_{V^*\times V}\dd r\\
 &= \frac{1}{2}\vert (S_h u')(s)\vert^2-\frac{1}{2}\vert (S_h u)'(t)\vert^2+ \calE^1(u_0+(KS_h u')(s))-\calE^1(u_0+(KS_h u')(t))\\
 &\quad +\calE^2_s(u_0+(KS_h u')(s))-\calE^2_t(u_0+(KS_h u')(t))
\end{align*} for all $s,t\in[0,T]$, where we have applied the integration by parts formula. This is possible as the duality pairing can be split now due to the fact that $(S_h u')(t)=\frac{1}{2h}(\tilde{u}(t+h)-\tilde{u}(t-h))$, where $\tilde{u}$ is a continuous extension of $u$ outside $[0,T]$. The continuous extension of $u$ is well-defined as $u\in \rmL^\infty(0,T;U)\cap \rmW^{1,1}(0,T;H)\subset \rmC_w([0,T];U)$ and therefore $S_h u'\in \rmL^2(0,T;U)$. However, we are not allowed to perform the limit passage after splitting up all the integrals: the duality pairing in the limit would not be well defined as we only know that $u''+\rmD\calE_t(u)\in \rmL^{q^*}(0,T;V^*)$. Nevertheless, since we have assumed $V\hookrightarrow \WW$, we can treat the term involving $\rmD\calE_t^2:\WW \rightarrow \WW^*\hookrightarrow V^*$ separately. First, taking into account 
\begin{align*}
u_0+(KS_h u')(t)=u_0+\frac{1}{2h}\int_{t-h}^{t+h} \tilde{u}(r)\dd r-\frac{1}{2h}\int_{-h}^{+h} \tilde{u}(r)\dd r
\end{align*} and that $u\in \rmC_w([0,T];U)\subset \rmC([0,T];\WW)$ since $U\overset{c}{\hookrightarrow} \WW$, there holds 
\begin{align}\label{eq:conv2}
\lim_{h\rightarrow 0}\left(u_0+(KS_h u') \right)=u \quad \text{in } \rmC([0,T];\WW). 
\end{align} Finally, by the continuity of $\calE_t^2$ and $\rmD\calE_t^2$, the convergences \eqref{eq:conv2} and $S_h u'\rightarrow u'$ in $\rmL^q(0,T;V)$ as $h\rightarrow 0$, there holds
\begin{align}\label{equality}
&=-\int_s^{t}\langle \rmD\calE^2_r(u(r)), u'(r)\rangle_{V^*\times V}\dd r\notag\\
&=\lim_{h\rightarrow 0}-\int_s^{t}\langle \rmD\calE^2_r(u_0+(KS_h u')(r)),(S_h u')(r)\rangle_{V^*\times V}\dd r\notag\\
 &= \lim_{h\rightarrow 0}\left(\calE^2_s(u_0+(KS_h u')(s))-\calE^2_t(u_0+(KS_h u')(t))\right)\notag\\
 &=\calE^2_s(u(s))-\calE^2_t(u(t))
\end{align} for all $s,t\in [0,T]$. Second, it has been shown in \cite[Lemma 6]{EmmTha11DNEE} that 
\begin{align*}
&-\int_0^{t}\langle u''(r)+\rmE(u(r)),u'(r)\rangle_{V^*\times V}\dd r\\
 &\leq \frac{1}{2}\vert v_0\vert^2-\frac{1}{2}\vert u'(t)\vert^2+ \calE^1(u_0)-\calE^1(u(t))
\end{align*} for almost every $t\in (0,T)$. The latter inequality  together with \eqref{equality} implies \eqref{formulaxy}, which completes the proof.

\end{proof}
 
\subsection{Proof of Theorem 2.4}\label{se:Proof.2}
\begin{proof}[\unskip\nopunct]
Let $u_0\in U, v_0\in H$ and $(\tau_n)_{n\in \mathbb{N}}$ be a vanishing sequence of positive step sizes. Let $(u_0^{k})_{k\in \mathbb{N}}\subset U\cap V$ and $(v_0^{k})_{k\in \mathbb{N}}\subset V$ be such that  $u_0^{k}\rightarrow u_0$ in $U$ and $v_0^{k}\rightarrow v_0$ in $H$ as $k\rightarrow \infty$. We let $k\in \mathbb{N}$ be fixed and we denote the interpolations associated with the initial data $u_0^k$ and $v_0^{k}$ again by \eqref{eq:Approx.U}-\eqref{eq:Approx.xif2} and suppress henceforth the dependence on $k$ for notational convenience. By Lemma \ref{le:LimitPass2}, there exists a subsequence (relabelled as before) of the interpolation functions and limit functions $u\in \rmC_w([0,T];U)\cap \rmW^{1,\infty}(0,T;H)\cap\rmW^{1,q}(0,T;V^*)\cap \rmW^{2,q^*}(0,T;U^*+V^*)$ (notice that $u_0^k \in U\cap V$) and $u(0)=u_0^k$ in $U$ and $u'(0)=v_0^{k}$ in $H$ such that the convergences \eqref{eq:LP2} hold, where we again suppress the dependence of the limit functions on $k$. First, we prove that the inclusion \eqref{sol:IC2} holds. To do so, we note that the \textsc{Euler--Lagrange} equation \eqref{eq:dis.inc2} reads 
\begin{align}\label{EuLa2}
\begin{split}
&\Vhn'(t)+\eta_{\tau_n}(t)+\rmD\calE_{\teon(t)}(\Uon(t))+S_{\ton}(t)=0 \quad \text{in $U^*+V^*$},\\
&\eta_n(t)\in \partial \Psi_{\Uun(t)}(\Von(t))
\end{split}
\end{align} for all $t\in (0,T)$, where $S_{\ton}(t)=B(\teon(t), \Vun(t),\Uun(t))-f_{\tau_n}(t), t\in [0,T]$. Testing equation \eqref{EuLa2} with $w\in \rmL^{\max\lbrace 2,q\rbrace}(0,T;U\cap V)$, we obtain
\begin{align*}
\int_0^T \langle \Vhn'(r)+\eta_{\tau_n}(r)+\rmD\calE_{\teon(s)}(\Uon(s))+S_{\ton}(r),w(r) \rangle_{(U^*+V^*)\times (U\cap V)}  \dd r=0.
\end{align*} Then, with the aid of the convergences \eqref{eq:LP2}, we are allowed to pass to the limit in the weak formulation obtaining
\begin{align*}
\int_0^T \langle u''(r)+\eta(r)+\rmD\calE_{s}(u(s))+B(t,u(r),u'(r))-f(r), w(r) \rangle_{(U^*+V^*)\times (U\cap V)} \dd r=0
\end{align*} for all $w\in \rmL^{\max\lbrace 2,q\rbrace}(0,T;U\cap V)$. Then, by a density argument and the fundamental lemma of calculus of variations, we deduce
\begin{align*}
u''(t)+\eta(t)+\rmD\calE_t(u(t))+B(t,u(t),u'(t))=f(t) \quad \text{in } U^*+V^*
\end{align*} for a.e. $t\in(0,T)$. We shall identify the weak limit $\eta$ as the subgradient of the dissipation potential, i.e, $\eta(t)\in \partial \Psi_{u(t)}(u'(t))$ for almost every $t\in (0,T)$. For that purpose, we will employ Lemma \ref{le:A1} with $f_n(t,v)=\Psi_{\Uun(t)}(v)$ and $f(t,v)=\Psi_{u(t)}(v)$ for all $v\in X=V$ and $n\in \mathbb{N}$. Assumption \eqref{le:A1.assum.1} is already satisfied by Condition \ref{eq:Psi2.4}. Hence, it remains to show 
\begin{align}\label{ineq:limsup}
\limsup_{n\rightarrow \infty}\int_0^T \langle \eta_n(t),\Von(t)\rangle_{V^*\times V}\dd t\leq \int_0^T \langle \eta(t),u'(t)\rangle_{V^*\times V}\dd t.
\end{align} In order to show \eqref{ineq:limsup}, we use the fact that $\eta_{\tau_n}$ satisfies the \textsc{Euler--Lagrange} equation \eqref{EuLa2}. Therefore, we will split the integral on the left-hand side of \eqref{ineq:limsup} and calculate 
\begin{align*}
&-\int_0^t \langle \Vhn'(r), \Von(r)\rangle_{V^*\times V} \dd r\\
&= -\int_0^t \langle \Vhn'(r), \Vhn(r)\rangle_{V^*\times V}\dd r+\int_0^t \langle \Vhn'(r), \Vhn(r)-\Von(r)\rangle_{V^*\times V}\dd r\\
&= \frac{1}{2}\vert v_0\vert^2-\frac{1}{2}\vert \Vhn(t)\vert^2+\int_0^t \langle \Vhn'(r), \Vhn(r)-\Von(r)\rangle_{V^*\times V}\dd r\\
&\leq  \frac{1}{2}\vert v_0\vert^2-\frac{1}{2}\vert \Vhn(t)\vert^2,
\end{align*} where we used the fundamental theorem of calculus for 
the absolutely continuous function $[0,T]\ni t\mapsto \frac{1}{2}\vert \Vhn(t)\vert^2$ and the estimate
\begin{align*}
&\int_0^t \langle \Vhn'(r), \Vhn(r)-\Von(r)\rangle_{V^*\times V}\dd r\\
&=\sum_{i=1}^{m-1} \int_{t_{i-1}}^{t_{i}}\left(\frac{V_{\tau_{n}}^i-V_{\tau_{n}}^{i-1}}{\tau_n},V_{\tau_{n}}^i\frac{r-t_{i-1}}{\tau_n}+V_{\tau_{n}}^{i-1}\frac{t_{i}-r}{\tau_n}-V_{\tau_{n}}^i\right)\dd r\\
&\quad+ \int_{t_{m-1}}^{t}\left(\frac{V_{\tau_{n}}^m-V_{\tau_{n}}^{m-1}}{\tau_n},V_{\tau_{n}}^m\frac{r-t_{m-1}}{\tau_n}+V_{\tau_{n}}^{m-1}\frac{t_{m}-r}{\tau_n}-V_{\tau_{n}}^m\right)\dd r\\
&=-\sum_{i=1}^{m-1} \int_{t_{i-1}}^{t_{i}}\left(\frac{V_{\tau_{n}}^i-V_{\tau_{n}}^{i-1}}{\tau_n},\left(V_{\tau_{n}}^i-V_{\tau_{n}}^{i-1}\right)\frac{t_{i}-r}{\tau_n}\right)\dd r \\
&\quad-\int_{t_{m-1}}^{t}\left(\frac{V_{\tau_{n}}^m-V_{\tau_{n}}^{m-1}}{\tau_n},\left(V_{\tau_{n}}^m-V_{\tau_{n}}^{m-1}\right)\frac{t_{m}-r}{\tau_n}\right)\dd r\\
&=-\sum_{i=1}^{m-1} \int_{t_{i-1}}^{t_{i}}\frac{t_{i}-r}{\tau_n^2}\left\vert V_{\tau_{n}}^i-V_{\tau_{n}}^{i-1}\right\vert^2\dd r -\int_{t_{m-1}}^{t}\frac{t_{m}-r}{\tau_n^2}\left\vert V_{\tau_{n}}^m-V_{\tau_{n}}^{m-1}\right\vert^2\dd r\leq 0
\end{align*} 
 for all $t\in (t_{m-1},t_m]$ and $m\in \lbrace 1,\dots,N\rbrace$.

We continue with the term involving the derivative of the energy functional and start with the linear part: 

\begin{align}\label{eq:est.E1}
&-\int_0^{t}\langle \rmE \Uon(r), \Von(r)\rangle_{U^*\times U}\dd r \notag \\
&= -\int_0^{t}\langle \rmE \Uhn(r), \Von(r)\rangle_{U^*\times U}\dd r+\int_0^{t}\langle \rmE \Uhn(r)-\rmE \Uon(r), \Von(r)\rangle_{U^*\times U}\dd r \notag\\
&= \calE^1(u_0)-\calE^1(\Uhn(t))+\int_0^{t}\langle \rmE (\Uhn(r)-\Uon(r)), \Von(r)\rangle_{U^*\times U}\dd r \notag\\
&\leq \calE^1(u_0)-\calE^1(\Uhn(t)),
\end{align} where we used
\begin{align*}
\int_0^{t}\langle \rmE (\Uhn(r)-\Uon(r)), \Von(r)\rangle_{U^*\times U}\dd r\leq 0,
\end{align*} which can be shown in the same way as above by using the strong positivity of $E$.
As for the nonlinear part, we obtain the desired estimate by employing the $\lambda$-convexity of $\calE^2_t$:
\begin{align}\label{eq:est.E2}
&-\int_0^{t}\langle \rmD\calE^2_{\teon(r)}(\Uon(r)), \Von(r)\rangle_{U^*\times U}\dd r \notag\\
&=-\sum_{i=1}^{m-1} \langle \rmD\calE^2_{t_i}(U_{\tau_n}^i), U_{\tau_n}^i-U_{\tau_{n}}^{i-1}\rangle_{U^*\times U}-\frac{t-t_{m-1}}{\tau_n} \langle \rmD\calE^2_{t_m}(U_{\tau_n}^m), U_{\tau_n}^m-U_{\tau_{n}}^{m-1}\rangle_{U^*\times U} \notag\\
&\leq -\sum_{i=1}^{m-1} \left(\calE^2_{t_{i}}(U_{\tau_n}^{i-1})-\calE^2_{t_i}( U_{\tau_n}^i)-\lambda\left\vert U_{\tau_n}^i-U_{\tau_{n}}^{i-1}\right\vert^2\right) \notag \\
&\quad -\frac{t-t_{m-1}}{\tau_n}  \left(\calE^2_{t_{m}}(U_{\tau_n}^{m-1})-\calE^2_{t_m}( U_{\tau_n}^m)-\lambda\left\vert U_{\tau_n}^m-U_{\tau_{n}}^{m-1}\right\vert^2\right) \notag\\
&= -\sum_{i=1}^{m} \left(\calE^2_{t_{i-1}}(U_{\tau_n}^{i-1})-\calE_{t_{i}}( U_{\tau_n}^i)+\int_{t_{i-1}}^{t_{i}}\partial_r\calE^2_r(U_{\tau_n}^i)\dd r+\lambda \tau_n^2\left\vert V_{\tau_n}^i \right\vert^2\right) \notag\\
&\quad +\frac{t_{m}-t}{\tau_n}  \left(\calE^2_{t_{m}}(U_{\tau_n}^{m-1})-\calE^2_{t_m}( U_{\tau_n}^m)-\lambda\left\vert U_{\tau_n}^m-U_{\tau_{n}}^{m-1}\right\vert^2\right) \notag \\
&= \calE^2_{0}(u_0)-\calE^2_{\teon(t)}(\Uon(t))+\int_{0}^{\teon(t)}\partial_r\calE^2_r(\Uon(r))\dd r+I_n(t),
\end{align} where 
\begin{align*}
I_n(t)=\frac{t_{m}-t}{\tau_n}  \left(\calE^2_{t_{m}}(U_{\tau_n}^{m-1})-\calE^2_{t_m}( U_{\tau_n}^m)-\lambda\left\vert U_{\tau_n}^m-U_{\tau_{n}}^{m-1}\right\vert^2\right)+\lambda \tau_n \int_{0}^{\teon(t)}\vert \Von(r)\vert^2\dd r.
\end{align*} Later, we will show that the term $I_n$ converges to zero as the step size vanishes. Next, we want to make use of the inequality \eqref{formulaxy}. However, the aforementioned inequality only holds true for almost every $t\in (0,T)$. Therefore, we take a sequence of increasing values $(\beta_l)_{l\in \mathbb{N}} \subset (0,T)$ that converges to $T$ for which the inequality \eqref{formulaxy} holds true.
Then, choosing $t=\beta_l$, we obtain with the convergences \eqref{eq:uhu.weak2}, \eqref{eq:vh.ptw.weak2}, \eqref{eq:B.uhh.ptw2}, and \eqref{eq:LP.vo.strong2}, the sequential weak lower semicontinuity of $\calE^1_t$ and $\vert\cdot\vert$, the continuity of $\calE^1_t$, Condition \ref{eq:cond.E2.4}, and \textsc{Fatou}'s Lemma that
\begin{align*} 
&\limsup_{n\rightarrow \infty} -\int_0^{\beta_l} \langle \Vhn'(r)+\rmD\calE_{\teon(r)}(\Uon(r)), \Von(r)\rangle_{V^*\times V} \dd r \notag \\
&\leq  \limsup_{n\rightarrow \infty}\Big( \frac{1}{2}\vert v_0 \vert^2-\frac{1}{2}\vert \Vhn(\beta_l) \vert^2 +\calE_{0}(u_0)-\calE^1(\Uhn(\beta_l))-\calE^2_{\teon(\beta_l)}(\Uon(\beta_l))\notag\\
&\quad +\int_{0}^{\teon(\beta_l)}\partial_r\calE_r(\Uon(r))\dd r+I_n(t) \Big)\notag \\
&\leq \frac{1}{2}\vert v_0 \vert^2-\frac{1}{2}\vert u'(\beta_l) \vert^2 + \calE_{0}(u_0)-\calE_{\beta_l}(u(\beta_l))+\int_{0}^{\beta_l}\partial_r\calE_r(u(r))\dd r,
\end{align*} where we also employed the inequalities \eqref{eq:est.E1} and \eqref{eq:est.E2} for the integral involving $\rmD \calE_{\teon(r)} = \rmD \calE^1+\rmD \calE_{\teon(r)}^2$. Taking into account that $u \in \rmC_w([0,T];U)$ and $u'\in \rmL^\infty(0,T;H)\cap \rmW^{1,1}(0,T;U^*+V^*)\subset\rmC_w([0,T];H)$, Lemma \ref{le:LimitPass2} yields
\begin{align*}
&\frac{1}{2}\vert v_0 \vert^2-\frac{1}{2}\vert u'(\beta_l) \vert^2+\calE_{0}(u_0)-\calE_{\beta_l}(u(\beta_l))+\int_{0}^{\beta_l}\partial_r\calE_r(u(r))\dd r\\
&\leq -\int_0^{\beta_l}\langle u''(r)+\rmD\calE_r(u(r)),u'(r)\rangle_{V^*\times V}.
\end{align*}  Then, in view of the convergences \eqref{eq:LP.f2} and \eqref{eq:LP.B2} as well as the \textsc{Euler--Lagrange} equation \eqref{EuLa2}, we obtain
\begin{align*}
&\limsup_{n\rightarrow \infty}\int_0^{\beta_l} \langle \eta_n(t), \Von(t)\rangle_{V^*\times V} \dd t\notag\\
&=\limsup_{n\rightarrow \infty}\int_0^{\beta_l} \langle S_{\tau_n}(t)-\Vhn'(t)-\rmD\calE_{\teon(t)}(\Uon(t)), \Von(t)\rangle_{V^*\times V} \dd t \notag\\
&\leq \int_0^{\beta_l} \langle f(t)-B(t,u(t),u'(t))-u''(t)-\rmD\calE_t(u(t)), u'(t) \rangle_{V^*\times V} \dd t \notag\\
&=\int_0^{\beta_l} \langle \eta(t), u'(t) \rangle_{U^*\times U} \dd t.
\end{align*} Together with Condition \ref{eq:Psi2.4} and Lemma \ref{le:A1}, this implies $\eta(t)\in \partial \Psi_{u(t)}(u'(t)) $ for almost every $t\in (0,\beta_l)$ for all $l\in \mathbb{N}$. Letting $l\rightarrow \infty$ leads to $\eta(t)\in \partial \Psi_{u(t)}(u'(t)) $ for almost every $t\in (0,T)$. For each $k\in \mathbb{N}$, this shows the existence of a function $u$ satisfying the inclusion \eqref{sol:IC2}, and the initial values $u(0)=u_0^k\in U\cap V$ and $u'(0)=v_0^{k}\in V$. \\ Now, denote with $(u_k)_{k\in \mathbb{N}}$ the sequence of solutions associated with sequence of initial values $(u_0^k)_{k\in \mathbb{N}}$ and $(v_0^k)_{k\in \mathbb{N}}$, and with $(\eta_k)_{k\in \mathbb{N}}$ the sequence of subgradients of $\Psi_{u_k(t)}(u_k'(t))$. In the last step, we want to show that there exists a limit function $u$ which satisfies \eqref{sol:IC2} and \eqref{sol:EDI2} as well as the initial values $u(0)=u_0$ in $U$ and $u'(0)=v_0$ in $H$. We recall that $u_0^{k}\rightarrow u_0$ in $U$ and $v_0^{k}\rightarrow v_0$ in $H$ as $k\rightarrow \infty$. 
We proceed the proof with the following steps:
\begin{itemize}
\item[1.] Deriving a priori estimates based on the energy-dissipation inequality \eqref{sol:EDI2},
\item[2.] Showing compactness of the sequences $(u_k)_{k\in \mathbb{N}}$ and $(\eta_k)_{k\in \mathbb{N}}$ in appropriate spaces,
\item[3.] Passing to the limit in the inclusion \eqref{sol:IC2} and the energy-dissipation balance \eqref{sol:EDI2} as $k\rightarrow \infty$.
\end{itemize}
\textbf{Ad 1.} 
Then, employing the convergences \eqref{eq:LP2}, we obtain
\begin{align*}
&\frac{1}{2}\vert u_k'(t) \vert^2 +\calE_t(u_k(t)) + \int_0^t \left(\Psi_{u_k(r)}(u_k'(r))+\Psi_{u_k(r)}^*(S_k(r)-\rmD \calE_{r}(u_k(r))-u_k''(r))\right)\dd r \\
&\leq \liminf_{n\rightarrow \infty}\bigg (
 \frac{1}{2} \left \vert \Von(t)\right \vert^2+ \calE_{\teon(t)}(\overline{U}_\ton(t)) \\
 &\quad  +\int_{0}^{\teon(t)}\left( \Psi_{\underline{U}_\ton(r)} ( \Von (r) ) + \Psi_{\underline{U}_\ton(t)}^*\left(S_\ton(r)- \rmD \calE_{\teon(r)}(\overline{U}_\ton(r))-\Vhn'(r)\right) \right) \dd r \bigg)
   \notag \\ 
  &\leq \limsup_{n\rightarrow \infty}\bigg ( \frac{1}{2} \left \vert v^k_0 \right \vert^2+
  \calE_{0}(u_0^k)+\int_{0}^{\teon(t)} \partial_r
  \calE_r(\Uun (r))\dd r\\
  &\quad+\int_{0}^{\teon(t)}
  \langle S_\ton(r), \Von(r) \rangle_{V^*\times V}\dd r
  +\tau \lambda\int_{0}^{\teo_\tau(t)}\Vert \Vo(r) \Vert_V^2\dd r\bigg ) \notag\\
 &= \frac{1}{2}\vert v_0^k \vert^2 +\calE_0(u_0^k) + \int_0^t  \partial_r \calE_r(u_k(r)) \dd r+\int_0^t \langle S_k(r),u_k'(r) \rangle_{V^*\times V}\dd r
\end{align*} for all $t\in [0,T]$, where $S_k(r)=f(r)-B(r,u_k(r),u_k'(r))$. Again, taking into account Conditions \ref{eq:cond.E2.4}, \ref{eq:B2.2}, and \ref{eq:f2}, we obtain with the classical lemma of \textsc{Gronwall}
\begin{align*}
&\frac{1}{2}\vert u_{k}'(t) \vert^2 +\calE_t(u_{k}(t)) + \int_0^t \left(\Psi_{u_k(r)}(u_{k}'(r))+\Psi_{u_k(r)}^*(S_{k}(r)-\rmD \calE_{r}(u_k(r))-u_{k}''(r))\right)\dd r \leq C.
\end{align*} for all $t\in [0,T]$ for a constant $C>0$.\\
\textbf{Ad 2.} With the same reasoning as for the interpolations, we obtain the convergences
\begin{subequations}
\label{eq:LP22}
\begin{align}
\label{eq:LP.uhu.weak.22}
u_{k} \overset{*}{\rightharpoonup} u \quad \text{in } &\rmL^{\infty}(0,T;U),\\
\label{eq:LP.u.22}
u_{k}-u_0^k \overset{*}{\rightharpoonup} u-u_0 \quad \text{in } &\rmL^{\infty}(0,T;V),\\
\label{eq:uhu.weak22}
u_k(t)\rightharpoonup u(t) \quad \text{in } &U \,\, \text{ for all }t\in[0,T],\\
\label{eq:uhu.weakinV22}
u_k(t)-u_0^k  \rightharpoonup u(t)-u_0 \quad \text{in } &V \,\, \text{ for all }t\in[0,T],\\
\label{eq:B.uhh22}
u_k \rightarrow u \quad \text{in } &\rmL^r(0,T;\WW) \quad \text{for any } r\geq 1,\\
\label{eq:B.uhh.ptw22}
u_k(t) \rightarrow u(t) \quad \text{in } &\WW \,\, \text{ for all }t\in[0,T],\\
\label{eq:LP.vo22}
u_k'(t) \overset{*}{\rightharpoonup} u' \quad \text{in } &\rmL^q(0,T;V)\cap \rmL^{\infty}(0,T;H),\\
\label{eq:LP.vo.strong22}
u'_k(t)\rightarrow u' \quad \text{in } &\rmL^p(0,T;H) \quad \text{for all }p\geq 1,\\
\label{eq:vh.ptw.weak22}
u_k'(t) \rightharpoonup u'(t) \quad \text{in } & H \,\, \text{ for all }t\in[0,T],\\
\label{eq:LP.eta2}
\eta^k_{\tau_n} \rightharpoonup \eta \quad \text{in } &\rmL^{q^*}(0,T;V^*),\\
\label{eq:LP.xi2.k}
\rmE u_k \rightharpoonup \rmE u \quad \text{in } &\rmL^2(0,T;U^*),\\
\label{eq:LP.xi22.k}
\rmD\calE^2_{t}(u_k) \rightarrow \rmD\calE^2_t(u) \quad \text{in } &\rmL^r(0,T;U^*)\quad \text{for any } r\geq 1,\\
\label{eq:LP.vhn'22}
u_k'' \rightharpoonup u'' \quad \text{in } &\rmL^{\min\lbrace 2,q^*\rbrace}(0,T;U^*+V^*),
\\
\label{eq:LP.B22}
B(\cdot,u_k,u_k')\rightarrow B(\cdot,u,u')\quad \text{in } &\rmL^{r^*}(0,T;V^*),
\end{align} 
\end{subequations} and if $\Psi_u$ satisfies \ref{eq:Psi2.3}, then 
\begin{subequations}
\begin{align*}
u'_k \rightarrow u' \quad \text{in } &\rmL^{\max\lbrace 2,q\rbrace}(0,T;V),\\
u_k \rightarrow u \quad \text{in } &\rmC([0,T];V).
\end{align*} 
\end{subequations}
\textbf{Ad 3.} Therefore, $ u\in\rmC_w([0,T];U)\cap  \rmW^{1,\infty}([0,T];H)\cap \rmW^{2,q^*}(0,T;U^*+V^*)$ with $u-u_0\in \rmW^{1,q}(0,T;V)$ and $\eta \in \rmL^{q^*}(0,T;V^*)$, and $u$ satisfies the initial conditions $u(0)=u_0$ in $U$ and $u'(0)=v_0$ in $H$. Along the same lines as for the interpolations, we obtain with Condition \ref{eq:Psi2.4} and Lemma \ref{le:A1} that $u$ and $\eta$ satisfy the inclusion \eqref{sol:IC2} and $u(t)\in \partial \Psi_{u(t)}(u'(t))$ a.e. in $(0,T)$. It remains to show the energy-dissipation balance \eqref{sol:EDI2}. The inequality 
\begin{align*}
&\frac{1}{2}\vert u'(t) \vert^2 +\calE_t(u(t)) + \int_0^t \left(\Psi_{u(r)}(u'(r))+\Psi_{u(r)}^*(S(r)-\rmD \calE_{r}(u(r))-u''(r)) \right)\dd r \\
&\leq  \frac{1}{2}\vert v_0 \vert^2 +\calE_0(u_0) + \int_0^t  \partial_r \calE_r(u(r)) \dd r+\int_0^t \langle S(r),u'(r) \rangle_{V^*\times V}\dd r,
\end{align*} for all $t\in [0,T]$ with $S(r)=f(r)-B(r,u(r),u'(r))$ is obtained by passing with $k\rightarrow \infty$ to the limit while taking into account the convergences \eqref{eq:LP2}. Then, employing again \eqref{formulaxy} and the \textsc{Fenchel--Young} inequality, we obtain
\begin{align*}
&\int_0^t \left(\Psi_{u(r)}(u'(r))+\Psi_{u(r)}^*(S(r)-\rmD \calE_{r}(u(r))-u''(r)) \right)\dd r\\
&\leq \frac{1}{2}\vert v_0 \vert^2-\frac{1}{2}\vert u'(t) \vert^2 +\calE_0(u_0)-\calE_T(u(t))+\int_0^t \partial_r \calE_r(u(r)) \dd r \\&\quad+\int_0^t \langle S(r),u'(r) \rangle_{V^*\times V}\dd r \\
&\leq  \int_0^t \langle \rmD\calE_r(u(r))-u''(r),u'(r) \rangle_{V^*\times V}\dd r+\int_0^t \langle S(r),u'(r) \rangle_{V^*\times V}\dd r\\
&=\int_0^t \langle S(r)-\rmD\calE_r(u(r))-u''(r),u'(r) \rangle_{V^*\times V}\dd r\\
&\leq \int_0^t \left(\Psi_{u(r)}(u'(r))+\Psi_{u(r)}^*(S(r)-\rmD \calE_{r}(u(r))-u''(r)) \right)\dd r
\end{align*} for almost every $t\in (0,T)$. Now, if $V\hookrightarrow U$, then the inequality \eqref{formulaxy} indeed holds as equality for all $t\in [0,T]$ by the classical integration by parts formula. This shows \eqref{sol:EDI2}, and thus the completion of the proof.
\end{proof}

\begin{rem} \label{rem.b} If we take a closer look at the proof, we notice  that we are allowed to consider the case $b\equiv 0$ if the compact embedding $V\overset{c}{\hookrightarrow}\widetilde{W}$ holds. This will ensure that we are able to pass to the limit as $\tau \searrow 0$ in the nonlinearity $\rmD\calE_t^2$.
\end{rem}

\begin{rem} The proof of Theorem \ref{th:MainExist2} reveals that one can consider dissipation potentials that depend on a parameter $\varepsilon$. In this case, the Condition \ref{eq:Psi2.1} is assumed to hold for every $\varepsilon \geq 0$ while the Condition \ref{eq:Psi2.2} holds uniformly in $\varepsilon\geq 0$. Condition \ref{eq:Psi2.4} can either be replaced with the \textsc{Mosco}-convergence $\Psi_{u_n}^{\varepsilon_n} \Mto \Psi_u^0$ for every sequence $u_n\rightharpoonup u$ as $\varepsilon\searrow 0$, or with a more general liminf estimate \eqref{le:A1.assum.1}, see \cite{BachoEmmrichMielke2019}.
\end{rem}

\begin{rem} It is feasible to consider nonsmooth energy functionals $\calE_t^2$ by regularizing the energy functional using the $p$-\textsc{Moreau--Yosida} \cite{Bach23GMYR} regularization and then passing to the limit with the regularized functional and its \textsc{G\^{a}teaux} derivative as the regularization parameter vanishes. In that case, we have to impose a growth condition and a time independence of $\calE^2$. We refer to \cite{Bach25DNEI} where this case has been executed in the case where the principal part of the operator acting on $u'$ is linear. This will extend our result to evolution inclusion with two multi-valued operators. Since this can be done in the exact same way as in \cite{Bach25DNEI}, this will not be considered here.
\end{rem}

\begin{section}{Applications}\label{se:App}
In this section, we want to apply the abstract result developed in the previous sections to concrete examples. First, we discuss in detail some physically motivated examples to illustrate the strength of the theory.

\section{Visco-elasto-plastic model for martensitic phase  transformation in shape-memory alloys}\label{ex.martensitic} In this example, we consider equations that describe a solid-solid phase transition in shape-memory alloys driven by stored-energy and a dissipation mechanism.
As critically discussed in \cite{RajRou03EDSA}, a commonly used model describing this phenomenon pertains to the isothermal case which is given by
\begin{align}\label{eq:marten}
\rho\partial_{tt}\uu+\nu(-1)^n \Delta^n \partial_t\uu -\nabla\cdot(\mathbold{\sigma}(\nabla \uu))+\mu (-1)^m \Delta^m \uu = \mathbold{f},
\end{align} where $\mathbold{f} \in \rmL^2(0,T;\rmH^{-1,2}(\Omega)^d)$, $m,n\in \mathbb{N}$ and $\mu,\nu\geq 0 $ are non-negative real values. Here, $\rho\geq 0$ denotes the density of the body, $\uu:\Omega \times[0,T] \rightarrow \mathbb{R}^d$ the displacement of the body, which is related to the deformation $\yy$ by $\uu(\xx,\cdot)=\yy(\xx,\cdot)-\xx$ on a reference body configuration $\Omega$, and $\mathbold{\sigma}:\mathbb{R}^{d\times d}\rightarrow \mathbb{R}^{d \times d}$ the \textsc{Piola--Kirchhoff} stress tensor depending on the gradient $\nabla u$. The stress tensor $\ssigma$ is, in general, not monotone, and for hyperelastic materials given by the derivative of a potential  $\varphi:\mathbb{R}^{d\times d}\rightarrow \mathbb{R}$ describing the specific stored energy, i.e, $\ssigma=\varphi'$. This implies that the potential $\varphi$ is in general not quasiconvex.\footnote{See, e.g.,  \cite[Remark 6.5, p. 175]{Roub13NPDE}.} The contribution of $\mu(-1)^m\Delta^m \uu$ in these equations models a capillarity-like behaviour of the solid, and $\nu(-1)^n\Delta^n \partial_t \uu$ describes a higher order viscosity. According to the authors, experiments show that the hysteretic phenomenon in shape memory alloys are rate-independent. This implies that the equation \eqref{eq:marten} does not model plasticity effects appropriately. The authors in \cite{RajRou03EDSA} suggest to incorporate a correction term into the equations which describes plasticity effects of the body. More precisely, they introduce a dissipation function $\lambda$ that is nonnegative and homogeneous of degree one that captures this hysteretic response. The governing equations are then given by
\begin{align*}
\rho\partial_{tt}\uu+\nu(-1)^n \Delta^n \partial_t\uu -\nabla\cdot(\ssigma_p+\mathbold{\sigma}(\nabla \uu))+\mu (-1)^m \Delta^m \uu = \mathbold{f},\\
\ssigma_p\in \mathrm{Sgn}\left( \lambda'(\nabla \uu): \nabla\partial_t\uu \right) \lambda'(\nabla \uu),
\end{align*} where $\mathrm{Sgn}: \mathbb{R}\rightrightarrows \mathbb{R}$ is  the multi-valued sign function, $\ssigma_p:\mathbb{R}^{d\times d}\rightarrow \mathbb{R}^{d\times d}$ is the plastic stress, and $\lambda:\mathbb{R}^{d\times d}\rightarrow \mathbb{R}$ is a so-called phase indicator and thus indicates the status of the phase of $\nabla \uu$. For a more physical discussion of the model, we refer to \cite{PleRou02VEPM,RajRou03EDSA,ArGrRo03MNSM} and the references therein where the model is studied extensively by \textsc{Roub{\'\i}{\v{c}}ek} and coauthors. In \cite{PleRou02VEPM}, the authors showed the existence of weak solutions for the case $n=m=2$ and $\nu, \mu >0$. In \cite{RajRou03EDSA}, the authors show the existence of very weak solutions for the critical cases $n=0$, $\nu,\mu \geq 0$ and $m\geq 3$ which can not be tackled in our framework due to the growth condition on $\Psi$ that requires $p>1$, see Condition \ref{eq:Psi2.4}.\\

With the theory developed here, we show the existence of solutions for the cases $\nu, \mu>0$ and $n\geq 1, m\geq 2$, which is not known in the literature under the assumptions on $\varphi$ presented here. We note that our theory also allows the case $\mu=0$ since $\nu>0$ creates enough regularity, see Remark \ref{rem.b}. For the sake of simplicity, we supplement the equation \eqref{eq:marten} with homogeneous \textsc{Dirichlet \& Neumann} boundary conditions. Non-homogeneity can be considered in a standard way, see \cite{PleRou02VEPM}. Hence, we consider the initial-boundary value problem 
\begin{align*}
\text{(P1)}
\begin{cases}
\rho\partial_{tt}\uu+\nu(-1)^n \Delta^n \partial_t\uu -\nabla\cdot(\ssigma_p+\mathbold{\sigma}(\nabla \uu))+\mu (-1)^m \Delta^m \uu = \mathbold{f} \quad \text{in } \Omega_T,\\
\ssigma_p(\xx,t)\in \mathrm{Sgn}\left( \lambda'(\nabla \uu(\xx,t)): \nabla\partial_t\uu(\xx,t) \right) \lambda'(\nabla \uu(\xx,t))\quad \text{a.e. in } \Omega_T,\\
\uu(\xx,0)\,\,=\uu_0(\xx) \quad \text{on } \Omega,\\
\uu'(\xx,0)=\vv_0(\xx) \quad \,\text{on } \Omega, \\
\frac{\partial^k \uu}{\partial \mathbold{\nu}^k}(\xx,t)\,\,=0 \quad  \qquad \text{on } \partial \Omega\times[0,T], \, k=0,\dots, \max\lbrace m,n\rbrace-1,
\end{cases}
\end{align*} where $\mathrm{Sgn}: \mathbb{R}\rightrightarrows \mathbb{R}$ is  the multi-valued and one-dimensional sign function, $\ssigma_p:\mathbb{R}^{d\times d}\rightarrow \mathbb{R}^{d\times d}$ is the plastic stress, and $\lambda:\mathbb{R}^{d\times d}\rightarrow \mathbb{R}$ is a so-called phase indicator and thus indicates the phase status of $\nabla \uu$. 
Moreover, $\rho:\mathbb{R}^d\rightarrow [0,\infty)$ is a  measurable function satisfying $\bar{\rho}\geq \rho(\xx)\geq\underline{\rho}>0$ for a.e. $\xx\in \Omega$. We want to show the existence of a weak solution to (P1) for any initial data $\uu_0\in \rmH_0^m(\Omega)^d$ and $\vv_0\in \rmL^2(\Omega)^d$ and external forces $\ff\in \rmL^2(0,T;\rmH^{-\max\lbrace m,n\rbrace}(\Omega)^d)$ in the following sense: there exists a function $\uu\in \rmC_w([0,T];\rmH_0^m(\Omega)^d)\cap \rmW^{1,\infty}(0,T;\rmL^2(\Omega)^d)\cap \rmH^{2}(0,T;\rmH^{-m}(\Omega)^d+\rmH^{-n}(\Omega)^d)$ with $u-u_0\in \rmH^{1}(0,T;\rmH^n_0(\Omega)^d)$ and $\ssigma_p\in \rmL^2(0,T;\rmH^{-n}(\Omega)^d))$ satisfying the initial conditions $\uu(0)=\uu_0$, $\vv(0)=\vv_0$, the integral equation
\begin{align}\label{weak.sense}
\begin{split}
&\int_0^T \left(  \langle \rho\uu'', \vv \rangle+\int_\Omega \left(\nu \nabla^n \partial_t\uu:\nabla^n \vv + (\ssigma_p+\ssigma(\nabla \uu)):\nabla \vv+ \mu \nabla^m \uu :\nabla^m \vv \right) \dd x\right)\dd t,\\
&=\int_0^T \langle \ff,\vv\rangle \dd t \quad \text{for all }\vv \in \rmL^2(0,T;\rmH_0^{\max\lbrace m,n\rbrace}(\Omega)^d),
\end{split}
\end{align} such that $\ssigma_p(\xx,t)\in \mathrm{Sgn}\left( \lambda'(\nabla \uu(\xx,t)): \nabla\partial_t\uu(\xx,t) \right) \lambda'(\nabla \uu(\xx,t))\quad \text{a.e. in } \Omega_T$, and the energy-dissipation balance 
\begin{align} \label{martensic.EDB}
& \frac{1}{2}\Vert \rho \uu'(t) \Vert_{\rmL^2(\Omega)^d}^2 +\frac{\mu}{2}\vert\uu(t)\vert_{m,2}^2 +\int_\Omega \varphi(\nabla \uu(t))\dd x +\int_\Omega\underset{r\in [0,t]}{\mathrm{Var}}(\lambda(\nabla \uu))\dd \xx+\int_0^t\frac{\nu}{2}\vert \uu'(r)\vert^2_{n,2}\dd r\notag \\&\quad+\int_0^t\Psi_{\uu(r)}^*(\ff(r)-\nabla\cdot(\ssigma_p+\ssigma (\nabla \uu(r)))+\mu (-1)^m\Delta^m \uu(r)-\uu''(r))\dd r \notag \\
&= \frac{1}{2}\Vert \rho \vv_0 \Vert_{\rmL^2(\Omega)^d}^2 + \frac{\mu}{2}\vert\uu_0\vert_{m,2}^2+\int_\Omega \varphi(\nabla \uu_0)\dd x+\int_0^t \langle \ff(t),\uu'(r) \rangle \dd r,
\end{align} for almost all $t\in(0,T)$, and if $n\geq m$, then the energy-dissipation balance holds for all $t\in(0,T)$. Here, $\Psi_{\uu}^*$ denotes the convex conjugate of $\Psi_{\uu}$ defined below. Moreover, $\underset{r\in [0,t]}{\mathrm{Var}}(\lambda)$ denotes the total variation of $\lambda$ over $[0,t]$ and $\langle \cdot,\cdot\rangle$ denotes the duality pairing between $\rmH_0^{\max\lbrace m,n\rbrace}(\Omega)^d$ and its dual space $\rmH^{-\max\lbrace m,n\rbrace}(\Omega)^d$, where $\rmH_0^{k}(\Omega)^d$ is the \textsc{Sobolev} space of all measurable functions whose weak derivatives exist up to the order $k\in \mathbb{N}$ and are square-integrable, and the traces of all derivatives up to the order $k-1$ vanish on the boundary $\partial \Omega$. It is readily seen that these spaces equipped with the inner product $(\vv,\ww)_{\rmH^k_0\times \rmH^{-k}}=\int_\Omega \nabla^k\vv:\nabla^k \ww \dd \xx$ form a \textsc{Hilbert} space, where $:$ is the \textsc{Frobenius} inner product.  It is well-known that by a classical density argument and the \textsc{Poincar\'{e}--Friedrichs} inequality, the norm induced by this inner product is equivalent to the standard norm. The (semi-)norm of $\rmH_0^{k}(\Omega)^d$ is denoted by $\vert \cdot\vert_{k,2}=\vert \cdot\vert_{\rmH^k_0(\Omega)}$ and with $\vert \cdot\vert_{-n,2}$ we denote the dual norm. Now, since the stored energy $\varphi$ was not supposed to satisfy any convexity assumption, we have in general two possibilities of approaching this problem. On the one hand, we can treat the stress tensor $\ssigma$ as strongly continuous perturbation of the capillarity if $\ssigma$ has at most linear growth. On the other hand, if we assume that the stress tensor satisfies an \textsc{Andrews--Ball} type condition allowing any polynomial growth for $\ssigma$, we can treat the stored energy $\varphi$ as part of the energy functional.\\\\
In this work, we choose the second approach and refer the interested reader to \cite{Bach21ONSA,PleRou02VEPM} for the first approach. Thus, we suppose that $\ssigma$ is a potential operator that satisfies an \textsc{Andrews--Ball} type condition. The \textsc{Andrews--Ball} type condition was originally introduced by \textsc{Andrews \& Ball} to show global existence of solutions for the one-dimensional equations in viscoelastodynamics, i.e., when $\nu>0, n=1$ and $\mu=0$, see  \cite{Andr80OESE,AndBal82ABCP}. The existence of weak solutions to the aforementioned case for arbitrary dimensions has already been studied in a more general abstract setting in \cite{EmmSis13EESO} by making the crucial assumption that the operator $B+\Lambda A$ is monotone for some $\Lambda>0$, which generalizes the \textsc{Andrews--Ball} condition. The \textsc{Andrews--Ball} condition  states that $\ssigma$ is monotone on a large scale, i.e., there exists a positive value $R>0$ such that
\begin{align*}
\left(\ssigma(\mathbold{F})-\ssigma(\mathbold{\tilde{F}})\right):\left( \mathbold{F}-\mathbold{\tilde{F}}\right)>0\quad\text{for all }\mathbold{F},\mathbold{F}\in \mathbb{R}^{d\times d} \text{ with }\vert \mathbold{F}-\mathbold{\tilde{F}}\vert\geq R,
\end{align*} where $\vert \cdot \vert$ is the norm induced by the \textsc{Frobenius} inner product.
We will impose the more general assumption of the convexity of $\varphi+\frac{\lambda}{2}$ which in the smooth setting is equivalent to the monotonicity of $\ssigma+\lambda \mathrm{id}$. However, if $m\in \mathbb{N}$ is sufficiently large so that we can again treat the stored energy $\varphi$ as strongly continuous perturbation, then the previous condition is redundant. Therefore, we will not explicitly focus on this case. Having said that, the exact conditions which we impose on the stress tensor $\ssigma$ and the phase indicator $\lambda$ are the following:
\begin{enumerate} [label=(\thesection.2\alph*), leftmargin=3.2em]  \label{sigma.II}
\item \label{sigma.b1} There exists a continuously differentiable function $\varphi:\mathbb{R}^{d\times d}\rightarrow \mathbb{R}$ such that $\ssigma=\varphi'$.
\item \label{sigma.b2} There exist positive constants $c^1_\sigma,C_\sigma^1>0$ and $p>1$ such that 
\begin{align*}
\vert \ssigma(\FF)\vert&\leq C_\sigma^1(1+\vert \FF\vert^{p-1})\\
c^1_\sigma\vert \FF \vert^{p}-C_\sigma^1\leq \vert \varphi(\FF)\vert&\leq C_\sigma^1(1+\vert \FF\vert^{p}) \quad \text{for all }\FF\in\mathbb{R}^{d\times d}.
\end{align*}
\item \label{sigma.b3} There exists a positive number $\Lambda>0$ such that $\varphi+\frac{\Lambda}{2}\vert\cdot\vert^2$ is convex. 
\item \label{lambda} There holds $\lambda\in \rmC^2(\mathbb{R}^{m\times d})$ such that $\lambda''$ is bounded.
\end{enumerate}

Condition \ref{sigma.b3} is in fact equivalent the $\Lambda$-convexity of $\varphi$ which in turn is equivalent to the following \textsc{Andrews--Ball} type condition: 
\begin{align}\label{andrew.balls}
\left(\ssigma(\mathbold{F})-\ssigma(\mathbold{\tilde{F}})\right):\left( \mathbold{F}-\mathbold{\tilde{F}}\right)\geq -\lambda \vert \mathbold{F}-\mathbold{\tilde{F}}\vert^2 \quad\text{for all }\mathbold{F},\mathbold{F}\in \mathbb{R}^{d\times d}.
\end{align} This follows from the convexity and \textsc{G\^{a}teaux} differentiability of $\varphi+\frac{\lambda}{2}\vert\cdot\vert^2$ and the parallelogram identity for $\vert\cdot\vert$. The \textsc{Andrews--Ball} condition in turn necessitates \eqref{andrew.balls} if $\ssigma$ is in addition locally \textsc{Lipschitz} continuous, see \cite{EmmSis13EESO}. Typically, $\varphi$ is a polynomial of order less than or equal to $4$ and is a multi-well potential, e.g., $\varphi(e)=(e^2-1)^2$ which is covered in our framework. As we mentioned before, this condition on the energy functional is more general than the ones considered in \cite{RajRou03EDSA} and \cite{PleRou02VEPM}, where the authors assume $\ssigma\in \rmC^2$ with $\ssigma''$ bounded and $\ssigma=\ssigma_1+\ssigma_2\in \rmC^1$ with $\ssigma_1$ being convex and $\ssigma_2$ being bounded, respectively. In particular, the case $\phi(e)=(e^2-1)^2=e^4+1-2e^2$ is not included in the aforementioned works. \\\\ The obvious choice for the function spaces are $U=\rmH_0^{m}(\Omega)^d,V=W=\rmH_0^{n}(\Omega)^d$, $\widetilde{W}=\rmW^{1,p}_0(\Omega)^d$ equipped with the standard norms, and $H=\rmL^2(\Omega)^d$ equipped with the inner product $(u,v)_{\rmL^2}=\int_\Omega \rho(x)u(x)v(x)\dd x$ and the induced norm. Furthermore, we assume $\ff\in \rmL^2(0,T;\rmH^{-\max\lbrace m,n\rbrace}(\Omega)^d)$. and choose $p> 1$ such that 
\begin{align}\label{embedding}
\rmH_0^{m}(\Omega)^d\overset{c}{\hookrightarrow}\rmW_0^{1,p}(\Omega)^d\hookrightarrow \rmL^{2}(\Omega)^d.
\end{align} This enables us to view the stored energy $\calE_2$ defined below as a strongly continuous perturbation of $\calE_1$. For example, we obtain for $m=2$ in dimension $d=2$ all values $p>1$ and in dimension $d=3$ the range $6/5\leq p\leq 6$. The dissipation potential $\Psi$ and the energy functional $\calE$ are given by
\begin{align*}
\Psi_{\uu}(\vv)=\int_{\Omega}\left(\frac{\nu}{2}\vert \nabla^n \vv(\xx)\vert^2 +\vert \lambda'(\nabla \uu(\xx)): \nabla \vv(\xx) \vert \right)\dd \xx=\Psi^1(\vv)+\Psi^2_{\uu}(\vv)
\end{align*} 
and
\begin{align*}
\calE(\uu)= \int_{\Omega} \left(\varphi(\nabla \uu(\xx))+\frac{\mu}{2}\vert \nabla^m \uu(\xx)\vert^{2}\right) \dd \xx=\calE_2(\uu)+\calE_1(\uu),
\end{align*} respectively, and therefore, $B\equiv 0$. Applying \cite[Corollary 3.5, p. 335]{EkeTem76CAVP} and \cite[Theorem 1.1, pp. 126]{AttBre86DSVF}, one can show that 
\begin{align*}
\Psi_{\uu}^*(\xii)=\min_{\mathbold{\eta}\in \rmH^1_0(\Omega)^d}\lbrace \min_{\overset{\mathbold{p}^*\in \rmL^2(\Omega)^d}{-\nabla\cdot \mathbold{p}^*=\mathbold{\eta}}} \int_\Omega f^*(\xx,\mathbold{p}^*(\xx))\dd \xx+\frac{1}{2 \nu}\vert \xii-\mathbold{\eta}\vert^2_{-n,2} \rbrace,
\end{align*} where $f^*$ is the real-valued convex conjugate of $f(\xx,\mathbold{p})=\vert \lambda'(\nabla \uu(\xx)):\mathbold{p}\vert$ with respect to $\mathbold{p}\in \mathbb{R}^{m\times d}$. Finally, we assume $n\geq 1,\nu,\mu >0$. \\\\
We start by verifying the assumptions on the dissipation potential $\Psi_u$.\\\\
\textbf{Assumptions on $\Psi_u$:} The dissipation potential $\Psi_u$ obviously complies with Condition \ref{eq:Psi2.1} as it is convex and finite everywhere on $V$. Since $\lambda'$ is supposed to be bounded, the growth condition \eqref{eq:Psi2.growth} in Condition \ref{eq:Psi2.2} is easily verified by the \textsc{Poincar\'{e}--Friedrichs} inequality. In order to prove Condition \ref{eq:Psi2.4}, we prove that for any sequence $u_n\rightharpoonup u$ as $n\rightarrow \infty$ with $\sup_{n\in \mathbb{N},t\in [0,T]}\calE_t(u_n)<+\infty$, there holds $\Psi_{u_n}\Mto \Psi_u$ is the sense of \textsc{Mosco}-convergence, see Remark \ref{re:Assump.Psi2} ii). As mentioned in the same remark, from \cite[Lemma 4.1]{Stef08BEPD}, we infer that $\Psi_{u_n}\Mto \Psi_u$ implies Condition \ref{eq:Psi2.4}. In order to show \eqref{Mosco}, we distinguish the cases $n=1$ and $n\geq 2$.\\\\
\textbf{Ad $n=1$:} Let $v_n\rightharpoonup v$ and $u_n\rightharpoonup u$ as $n\rightarrow \infty$ with $\sup_{n\in \mathbb{N},t\in [0,T]}\calE_t(u_n)<+\infty$. By the compact embedding \eqref{embedding}, there exists a subsequence (labelled as before) and a function $g\in \rmW^{1,p}(\Omega)^d$ such that 
\begin{align}\label{convergence.u}
\nabla u_n(\xx)\rightarrow \nabla u(\xx)\quad &\text{for a.e. }\xx\in \Omega,\\
\vert \nabla u_n(\xx)\vert\leq g(\xx)\quad &\text{for a.e. }\xx\in \Omega.
\end{align}

We define 
\begin{align*}
f(\xx,\zz,\xii)=f(\zz,\xii)=\vert \lambda'(\zz)\xii\vert +\vert \xii\vert^2
\end{align*} and note that since $\lambda'$ is continuous, $f$ satisfies the assumptions of Theorem 2.1 in \cite{EkeTem76CAVP}, which implies 
\begin{align}\label{conv.a}
\Psi_{u}(v)&=\int_\Omega \left( \vert \nabla v(\xx)\vert^2+\vert \lambda'(\nabla u(\xx)):\nabla v(\xx)\vert \right) \dd \xx \notag \\
&\leq \liminf_{n\rightarrow \infty} \int_\Omega \left( \vert \nabla v_n(\xx)\vert^2+\vert \lambda'(\nabla u_n(\xx)):\nabla v_n(\xx)\vert \right) \dd \xx \notag \\
&= \liminf_{n\rightarrow \infty} \Psi_{u_n}(v_n).
\end{align} From the fact that $\lambda''$ is bounded it follows that $\lambda'$ has at most linear growth. Hence, we obtain for the constant sequence $\tilde{v}_n=v,$ $n\in \mathbb{N}$, by the dominated convergence theorem 
\begin{align}\label{conv.b}
\Psi_{u}(v)&=\int_\Omega \left( \vert \nabla v(\xx)\vert^2+\vert \lambda'(\nabla u(\xx)):\nabla v(\xx)\vert \right) \dd \xx \notag \\
&= \lim_{n\rightarrow \infty} \int_\Omega \left( \vert \nabla v(\xx)\vert^2+\vert \lambda'(\nabla u_n(\xx)):\nabla v(\xx)\vert \right) \dd \xx \notag \\
&= \lim_{n\rightarrow \infty}\Psi_{u_n}(v)
\end{align} from which the \textsc{Mosco}-convergence \eqref{Mosco} follows. In a standard way, one shows the convergence of the whole sequence.\\\\
\textbf{Ad $n\geq 2$:} Let $v_n\rightharpoonup v$ in $V$ and $u_n\rightharpoonup u$ in $U$ as above. Since $n\geq 2$, by the compact embedding $\rmH_0^{n}(\Omega)^d\overset{c}{\hookrightarrow} \rmH_0^{1}(\Omega)^d$, there holds $v_n\rightarrow v$ in $\rmH^{1}_0(\Omega)^d$. Then, taking again \eqref{convergence.u} into account, we obtain \eqref{conv.a} and \eqref{conv.b} with the dominated convergence theorem, whence Condition \eqref{Mosco}.\\\\
\textbf{Assumptions on $\calE$:} 
Now, we wish to verify the Conditions \ref{eq:cond.E2.1}-\ref{eq:cond.E2.8} for $\calE$. From the assumptions \ref{sigma.b1}-\ref{lambda} and the fact that $\calE$ is time-independent, Conditions \ref{eq:cond.E2.1}-\ref{eq:cond.E2.4} are easily verified. From Assumption \ref{sigma.b3} and the parallelogram identity of the norm of $\rmH^m_0(\Omega)^d$, it follows that 
\begin{align*}
\calE(\vartheta u+(1-\vartheta )v)\leq &\vartheta \calE(u)+(1- \vartheta)\calE(v)+\vartheta(1-\vartheta)\left(\lambda\vert u-v\vert^2_{\rmH_0^1(\Omega)^d}-\mu \vert u-v\vert^2_{\rmH_0^m(\Omega)^d}\right)
\end{align*} for all $t\in[0,T],  \vartheta\in [0,1]$ and $u,v\in U$. Employing the \textsc{Gagliardo--Nirenberg} inequality, see, e.g., \cite{Nire66AEIE}, \cite{Frie69PDE} or \cite[Section 21.19]{Zeid90NFA2a}, there exist constants $c_1,c_2>0$ such that 
\begin{align*}
\vert u-v\vert^2_{\rmH_0^1(\Omega)^d}&\leq c_1\vert u-v\vert_{\rmH^m_0(\Omega)^d}^{2/m}\Vert u-v\Vert_{\rmL^2(\Omega)^d}^{2(m-1)/m}+c_2\Vert u-v\Vert_{\rmL^2(\Omega)^d}^2\\
&\leq \varepsilon \vert u-v\vert_{\rmH^m_0(\Omega)^d}^{2}+C(\varepsilon) \Vert u-v\Vert_{\rmL^2(\Omega)^d}^{2}+c_2\Vert u-v\Vert_{\rmL^2(\Omega)^d}^2,
\end{align*} where we employed \textsc{Young}'s inequality in the last step for $\varepsilon>0$ and a constant  $C(\varepsilon)>0$. Choosing a sufficiently small $\varepsilon>0$ (e.g. $\varepsilon<\mu$), we obtain the $\tilde{\Lambda}$-convexity of $\calE$ for $ \tilde{\Lambda}:=C(\varepsilon)+c_2$. Thus, the energy functional $\calE$ fulfills Condition \ref{eq:cond.E2.7}. Further,  Assumptions \ref{sigma.b1} - \ref{lambda} imply that $\calE_1$ and $\calE_2$ are \textsc{Fr\'{e}chet} differentiable on $U$ and $\widetilde{W}$, respectively, and the \textsc{Fr\'{e}chet} derivatives are given by
\begin{align*}
\langle \mathrm{D}\calE_2(\uu),\vv\rangle_{U^*\times U}=\mu\int_\Omega \nabla^m \uu(\xx)\cdot \nabla^m \vv(\xx)\dd \xx\\
\langle \mathrm{D}\calE_1(\uu),\vv\rangle_{\widetilde{W}^*\times \widetilde{W}}=\int_\Omega \ssigma(\nabla \uu(\xx)):\nabla \vv(\xx) \dd \xx.
\end{align*} Consequently, by the subdifferential calculus, the subdifferentials are single-valued with $\partial \calE_1(\uu)= \lbrace \mathrm{D}\calE_1(\uu)\rbrace$ and $\partial \calE_2(\uu)= \lbrace \mathrm{D}\calE_2(\uu)\rbrace$. Hence, $\xii_1 \in \partial \calE_1(\uu)$ and $\xii_2 \in \partial \calE_2(\uu)$ if and only if $\xii_1=\mu (-1)^m\Delta^m \uu \in U^*$ and $\xii_2=-\nabla\cdot\ssigma (\nabla \uu)\in \widetilde{W}^*$, respectively.
Then, Condition \ref{eq:cond.E2.8} follows from the following estimate:
\begin{align*}
\langle \mathrm{D}\calE_2(\uu),\vv\rangle_{{\widetilde{W}}^*\times {\widetilde{W}}}&=\int_\Omega \ssigma(\nabla \uu(\xx)):\nabla \vv(\xx)\dd \xx\\
&\leq \left(\int_\Omega \vert \ssigma(\nabla \uu(\xx))\vert^{p/(p-1)}\dd \xx\right)^{(p-1)/p}\left(\int_\Omega\vert\nabla \vv(\xx)\vert^p\dd \xx\right)^\frac{1}{p}\\
&\leq \left(C_\sigma 2^{p/(p-1)}\int_\Omega (1+\vert \nabla \uu(\xx)\vert^{p})\dd \xx\right)^{(p-1)/p}\Vert \vv\Vert_{\widetilde{W}}\\
&\leq C\left(1+ \int_\Omega \vert \nabla \uu(\xx)\vert^{p}\dd \xx\right)\Vert \vv\Vert_{\widetilde{W}}\\
&\leq C\left (1+\int_\Omega \varphi(\nabla \uu(\xx))\dd \xx\right)\Vert \vv\Vert_{\widetilde{W}}\\
&\leq C(1+\calE_2(\uu))\Vert \vv\Vert_{\widetilde{W}},
\end{align*} where we have employed \textsc{H\"{o}lder}'s and \textsc{Young}'s inequality as well as the growth condition \ref{sigma.b2}. Again, $C>0$ denotes a generic constant thah is independent of $u$ and that can change from line to line. This shows Condition \ref{eq:cond.E2.8}. By \textsc{H\"{o}lder}'s inequality, we also obtain 
\begin{align*}
\langle \mathrm{D}\calE_2(\uu)-\mathrm{D}\calE_2(\vv),\ww\rangle_{{\widetilde{W}}^*\times {\widetilde{W}}}&=\int_\Omega \left(\ssigma(\nabla \uu(\xx))-\ssigma(\nabla \vv(\xx))\right)\nabla \ww(\xx)\dd \xx\\
&\leq \left(\int_\Omega \vert \ssigma(\nabla \uu(\xx))-\ssigma(\nabla \vv(\xx))\vert^{p/(p-1)}\dd \xx\right)^{(p-1)/p}\Vert \ww \Vert_{\widetilde{W}}
\end{align*} and therefore
\begin{align*}
\Vert \mathrm{D}\calE_2(\uu)-\mathrm{D}\calE_2(\vv)\Vert_{\widetilde{W}^*} \leq \left(\int_\Omega \vert \ssigma(\nabla \uu(\xx))-\ssigma(\nabla \vv(\xx))\vert^{p/(p-1)}\dd \xx\right)^{(p-1)/p}.
\end{align*} Then, Condition \ref{eq:cond.E2.6} follows from the dominated convergence theorem.
Since all conditions of Theorem \ref{th:MainExist2} are verified, for every $\uu_0\in U$ and $\vv_0\in H$ there exists a weak solution $\uu\in \rmC_w([0,T];\rmH_0^m(\Omega)^d)\cap \rmW^{1,\infty}(0,T;\rmL^2(\Omega)^d)\cap \rmH^{2}(0,T;\rmH^{-m}(\Omega)^d+\rmH^{-n}(\Omega)^d)$ with $u-u_0\in \rmH^{1}(0,T;\rmH^n_0(\Omega)^d)$ and $\ssigma_p\in \rmL^2(0,T;\rmH^{-n}(\Omega)^d))$ satisfying the integral equation \eqref{weak.sense}. Noting that there holds
\begin{align*}
\underset{r\in [0,t]}{\mathrm{Var}}\lambda(\nabla \uu))=\int_0^t\left \vert \frac{\partial}{\partial r}\lambda(\nabla \uu(r))\right\vert\dd r=\int_0^t \left \vert \lambda'(\nabla \uu(r)):\nabla \partial_t\uu(r)\right\vert\dd r=\int_0^t\Psi^2_{\uu}(\uu'(r))\dd r,
\end{align*} we infer the energy-dissipation balance \eqref{martensic.EDB}.

\begin{rem} We note that this result has not been shown before in the literature and that there are no abstract results that can address this type of problem. 
\end{rem}

\subsection{Differential Inclusion II}\label{ex:Klein.Gordon} The following example is a nonlinearly damped inertial system and can, for smooth dissipation potentials, be interpreted as a viscous regularization of the \textsc{Klein--Gordon} equation. The equations supplemented with initial and boundary conditions read
\begin{align*}
\text{(P2)}
\begin{cases}
\partial_{tt} u-\nabla \cdot \mathbf{p} -\Delta u +b(u) = f \quad \text{in } \Omega_T,\\
\mathbf{p}(\xx,t)\in \partial_v \psi(\xx,u(\xx,t), \nabla \partial_t u(\xx,t)) \quad \text{a.e. in } \Omega_T,\\
u(\xx,0)\,\,=u_0(\xx) \quad \text{on } \Omega,\\
u'(\xx,0)=v_0(\xx) \quad \,\text{on } \Omega, \\
u(\xx,t)\,\,=0 \quad  \qquad \text{on } \partial \Omega\times[0,T].
\end{cases}
\end{align*} If $\psi=0$ and $b(u)=\gamma u$ for a constant $\gamma>0$, then the equation in (P2) reduces to the classical \textsc{Klein--Gordon} equation, which is a relativistic wave equation that is related to the \textsc{Schr\"{o}dinger} equation and has applications in relativistic quantum mechanics.\\

We make the following assumptions on the functions $\psi$ and $b$. For simplicity, we choose $d=1$ and note that the case $d\geq 2$  can, under stronger assumptions, be treated in a similar way.
\begin{enumerate} [label=(\thesection.\alph*), leftmargin=3.2em]  \label{eq:psi.applic.II}
\item \label{2.a} The function $\psi:\Omega \times \mathbb{R}\times \mathbb{R} \rightarrow [0,+\infty)$ is a \textsc{Carath\'{e}odory} function such that $\psi(x,y, \cdot)$ is a proper, lower semicontinuous, and convex functional, and there holds $\psi(y,y, 0)=0$ for almost every $x\in \Omega$ and all $y\in \mathbb{R}$.
\item \label{2.b} There exist a real number $q>1$ and positive constants $ c_\psi,C_\psi>0$ such that 
\begin{align*}
c^R_\psi\left(\vert z\vert^q-1\right)\leq \psi(x,y,z)&\leq C^R_\psi\left( 1+\vert z\vert^q\right) \quad 
\end{align*} $\text{for a.e. }x\in \Omega \text{ and all }z \in \mathbb{R}^m$ and  $y\in \mathbb{R}$ with $\vert y\vert \leq R$. 
\item \label{2.d} The function $b:\Omega\rightarrow \mathbb{R}$ is continuous and there exist a real number $p>1$ and a constant $C_b>0$ such that
\begin{align*}
\vert b(u)\vert\leq C_b(\vert u\vert^{p-1}+1) \quad \text{ for all }u\in \mathbb{R}.
\end{align*}
\end{enumerate}
Accordingly, the function spaces are given by $V=\rmW^{1,q}_0(\Omega),\, U=\rmH_0^1(\Omega),\WW=\rmL^{\max\lbrace p,2\rbrace}(\Omega)$ and $H=\rmL^2(\Omega)$. Then, we identify the dissipation potential $\Psi:V\rightarrow \mathbb{R}$ and the energy functional $\calE:U\rightarrow [0,+\infty)$ by
\begin{align*}
\Psi_u(v)=\int_{\Omega} \psi(x,u(x),\nabla v(x)) \dd x \quad \text{and} \quad \calE(u)=\frac{1}{2}\int_{\Omega} \vert \nabla u(x)\vert^{2} \dd x,
\end{align*} respectively. The perturbation $B:\WW \rightarrow V^*$ is given by
\begin{align*}
\langle B(u),w \rangle_{\WW^*\times \WW}= \int_{\Omega} b(u(x)) w(x) \dd x.
\end{align*} We note that the conjugate functional $\Psi_u^*$ can not, in general, be expressed as an integral functional over $\Omega$, since it is defined on $\rmW^{-1,q^*}(\Omega)$.\\

Obviously, $\calE$ satisfies all Conditions \ref{eq:cond.E2.1}-\ref{eq:cond.E2.8}. In view of the compact embedding $\rmH_0^1(\Omega)\overset{c}{\hookrightarrow}\rmC(\overline{\Omega})$ and \textsc{Fatou}'s lemma, it is readily seen that $\Psi_u$ satisfies Conditions \ref{eq:Psi2.1} and \ref{eq:Psi2.2}. In order to verify Condition \ref{eq:Psi2.4}, we show that for every sequence $u_n\rightharpoonup u$ in $U$ with $\sup_{n\in \mathbb{N}}\calE(u_n)<+\infty$, there holds $\Psi_{u_n}\Mto \Psi_u$ as $n\rightarrow \infty$. As we mentioned in Remark \ref{re:Assump.Psi2} ii), the \textsc{Mosco}-convergence  $\Psi_{u_n}\Mto \Psi_u$ implies Condition \ref{eq:Psi2.4}. The liminf estimate in the \textsc{Mosco}-convergence follows from \cite[Theorem 3]{Ioff77LSIF}. The limsup estimate is trivially fulfilled by choosing, for each $v\in V$, the constant sequence $v_n=v, n\in \mathbb{N}$, and the dominated convergence theorem.

If we assume $p\in (1,2]$, it is easy to check in the same way as in the previous example that Conditions \ref{eq:B2.1} and \ref{eq:B2.2} are also fulfilled. Finally, we assume $f\in \rmL^{2}(0,T;H)$, $u_0\in U$ and $v_0\in H$. Therefore, Theorem \ref{th:MainExist2} ensures the existence of a solution $u\in \rmC_w([0,T];U)\cap \rmW^{1,\infty}(0,T;H)\cap \rmW^{2,q^*}(0,T;U^*+V^*)$ with $u-u_0\in \rmW^{1,q}(0,T;V)$ to (P2), satisfying the integral equation
\begin{align*}
\int_0^T \Big( \langle u'' v\rangle_{(U^*+V^*)\times(U\cap V)} +\int_\Omega    \mathbf{p}\cdot \nabla v +b(u)v \dd x \Big) \dd t = \int_0^T\int_\Omega fv\dd x \dd t 
\end{align*} \text{for all }$v\in \rmL^{\min\lbrace 2,q^*\rbrace}(0,T;U^*+V^*)$ with $\mathbf{p}(x,t)\in \partial_v \psi(x,u(x,t), \nabla \partial_t u(x,t))$ a.e. in $\Omega_T$. Furthermore, the energy-dissipation balance 
\begin{align*}
&\frac{1}{2}\Vert u'(t)\Vert^2_{\rmL^2(\Omega)}+\frac{1}{2}\Vert u(t)\Vert^2_{H_0^1(\Omega)} +\int_0^t \left( \Psi_{u(t)}(u'(r))+\Psi_{u(t)}^*({f}(r)-u''(r)-\Delta u(r))\right) \dd r \notag \\
&= \frac{1}{2}\Vert v_0\Vert^2_{\rmL^2(\Omega)}+\frac{1}{2}\Vert u_0\Vert^2_{H_0^1(\Omega)} +\int_0^t \langle f(r), u'(r)\rangle_{\rmL^2(\Omega)\times\rmL^{2}(\Omega)} \dd r
\end{align*} holds for almost every $t\in (0,T)$ when $q\in(1,2)$ and for all $t\in (0,T)$ when $q\geq 2$.
\subsection{Differential inclusion III} \label{se:app.DII}
In the final example, we consider a nonlinearly damped inertial system which can be interpreted as a model in ferro-magnetism where inertia is taking into account \cite{MiRoSa13NADN}. The differential inclusion supplemented with initial and boundary conditions is given by 
\begin{align*}
\text{(P3)}
\begin{cases}
\partial_{tt} u+\left \vert \partial_t u\right \vert^{q-2}\partial_t  u +p -\nabla \cdot \left( E\nabla u\right) +W'(u) = f \quad \text{in } \Omega_T,\\
p(\xx,t)\in \mathrm{Sgn}\left(\partial_t u(\xx,t)\right)\quad\text{a.e. in } \Omega_T,\\
u(\xx,0)\,\,=u_0(\xx) \quad \text{on } \Omega,\\
u'(\xx,0)=v_0(\xx) \quad \,\text{on } \Omega, \\
u(\xx,t)\,\,=0 \quad  \qquad \text{on } \partial \Omega\times[0,T],
\end{cases}
\end{align*} where $q\geq 2$, $\rmW:\mathbb{R}\rightarrow \mathbb{R}$ is a $\lambda$-convex and continuously differentiable function, $\rmE:\mathbb{R}^m\rightarrow \mathbb{R}^m$ is a uniformly positive definite and symmetric matrix, and $\mathbold{f}\in \rmC^1([0,T];\rmH^{-1}(\Omega))$. As mentioned before, the double-well potential $W(\uu)=(1-\uu^2)^2$ is a admissible choise for $W$. We set $U=\rmH_0^1(\Omega),\, V=\rmL^{q}(\Omega),$ and $H=\rmL^2(\Omega)$. Then, the dissipation potential $\Psi:V\rightarrow \mathbb{R}$ and the energy functional $\calE:U\rightarrow [0,+\infty]$ are given by
\begin{align*}
\Psi_u(v)&=\Psi(v)=\int_{\Omega}\left(\frac{1}{q}\vert v(\xx)\vert^q+\vert v(\xx)\vert\right)\dd \xx \, \text{ and } \\ \calE_t(u)&=\int_{\Omega} \frac{1}{2}\nabla u(\xx):E(x)\nabla u(\xx) \dd \xx+\int_{\Omega} W(u(\xx)) \dd \xx-\langle f(t),u\rangle_{U^*\times U}
\end{align*} respectively. Consequently, $B=0$ and 
\begin{align*}
    \calE_t^2(u)= \int_{\Omega} W(u(\xx)) \dd \xx-\langle f(t),u\rangle_{U^*\times U}.
\end{align*} By our assumptions, it is easy to see hat the conditions of Theorem \ref{th:MainExist2} are satisfied. Therefore, for all $u_0\in U$ and $v_0\in H$, there exists a solution $u\in \rmC_w([0,T];U)\cap \rmW^{1,\infty}(0,T;H)\cap \rmW^{2,q^*}(0,T;U^*+V^*)$ with $u-u_0\in \rmW^{1,q}(0,T;V)$ to (P3), i.e., $u$ fulfills the integral equation
\begin{align*}
&\int_0^T \Big( \langle u'' v\rangle_{(U^*+V^*)\times(U\cap V)} +\int_\Omega  \vert \partial_t u\vert^{q-2} \partial_t u v +pv +\nabla u\cdot\nabla v \dd \xx \Big) \dd t =\int_0^T\langle f,u\rangle_{U^*\times U} \dd t 
\end{align*} for all $v\in \rmL^{\min{\lbrace 2, q^*\rbrace}}(0,T;U^*+V^*)$ with ${p}(t,\xx)\in \mathrm{Sgn}( u(\xx,t))$ a.e. in $\Omega_T$, and the energy-dissipation balance \eqref{sol:EDI2} holds for almost every $t\in (0,T)$.
\end{section}

\appendix
\section{Appendix}
\label{se:Appendix}

\subsection{Subdifferential calculus}

In this section, we want to collect some of the results from the theory of subdifferential calculus. Let $(X,\Vert \cdot\Vert)$ be a separable and reflexive \textsc{Banach} space and denote with $(X^*,\Vert\cdot\Vert_*)$ its topological dual space.  Unlike the differential operator, the subdifferential operator is, in general, not linear which often causes technical difficulties. The following well-known result shows under which assumptions on the functionals, linearity of the subdifferential operator holds, also known as the variational sum rule. 
\begin{lem}[Variational sum rule] \label{le:Subdif} 
\begin{itemize}
\item[1)] Let $f_1:X\rightarrow (-\infty,+\infty]$ and $f_2:X\rightarrow (-\infty,+\infty]$ be subdifferentiable and let $f_2$ be \textsc{Fr\'{e}chet} differentiable in $u\in \DOM(\partial f_1)\cap \DOM(\partial f_2)\neq \emptyset$. Then, 
\begin{align*}
\partial(f_1+f_2)(u)=\partial f_1(u)+D f_2(u),
\end{align*} where $D f_2$ denotes the \textsc{Fr\'{e}chet} derivative of $f_2$.
\item[2)] Let $f_1:X\rightarrow (-\infty,+\infty]$ and $f_2:X\rightarrow (-\infty,+\infty]$ be proper, lower semicontinuous and convex, and if there is a point $\tilde{u}\in \DOM(f_1)\cap \DOM(f_2)$ where $f_2$ is continuous, there holds 
\begin{align}
\partial (f_1+f_2)(v)=\partial f_1(v)+\partial f_2(v) \quad \text{for all }v\in X.
\end{align}
If $f_2$ is in addition \textsc{G\^{a}teaux} differentiable on $X$, there holds $\partial f_2(v)=\lbrace D_G f_2(v)\rbrace$ and we have
\begin{align*}
\partial(f_1+f_2)(v)=\partial f_1(v)+D_G f_2(v) \quad \text{for all } v\in X,
\end{align*} where $D_G f_2$ denotes the \textsc{G\^{a}teaux} derivative of $f_2$.
\end{itemize}
\end{lem}

\begin{proof}
The proof of assertion 1) follows immediately from the definition of a subdifferential and the proof of assertion 2) follows from Proposition 5.3. on p. 23 and Proposition 5.6 on p. 26 in \cite{EkeTem76CAVP}.
\end{proof}

The next lemma establishes a deep connection between the subgradient of a functional and its convex conjugate $f^*(\xi):=\sup_{u\in X}\left \lbrace \langle \xi,u\rangle -f(u)\right \rbrace, \quad \xi\in X^*$.
 
\begin{lem}\label{le:Leg.Fen}
 Let $X$ be a \textsc{Banach} space and let $f:X\rightarrow (-\infty,+\infty]$ be a proper, lower semicontinuous, and convex functional and let $f^*:X^*\rightarrow (-\infty,+\infty]$ be the convex conjugate of $f$. Then for all $(u,\xi)\in X\times X^*$, the following assertions are equivalent:
\begin{itemize}
\item[$i)$]$\xi\in \partial f(u) \quad \text{in } X^*;$
\item[$ii)$] $u\in \partial f^*(\xi)\quad \text{in } X;$
\item[$iii)$]$\langle \xi, u\rangle_{X^*\times X}=f(u)+f^*(\xi) \quad \text{in
  } \mathbb{ R}.$
\end{itemize}
\end{lem}
\begin{proof} Proposition 5.1 and Corollary 5.2 on pp. 21 in \cite{EkeTem76CAVP}.
\end{proof}
In the next result, we show under which conditions on a time-dependent functional $f:[0,T]\times X\rightarrow (-\infty,+\infty]$, the associated time-integral functional $F$ 
\begin{flalign} \label{ap:integral.functional}
F(x)=\begin{cases}
\int_0^T f(t,\xi(t))\dd t \quad &\text{if } f(\cdot,\xi(\cdot))\in \rmL^1(0,T),\\
+\infty \quad &\text{otherwise}.
\end{cases}
\end{flalign} inherits the properties of $f$. Before, we introduce some notions and definitions. We denote with $\mathscr{L}_{(0,T)}$ the \textsc{Lebesgue} $\sigma$-algebra of the interval $[0,T]$ and with $\mathscr{B}(X)$ the \textsc{Borel} $\sigma$-algebra of $X$. A functional $f:[0,T]\times
X\rightarrow (-\infty,+\infty]$ is called a \textit{normal integrand} if it is $\mathscr{L}_{(0,T)} \otimes
\mathscr{B}(X)$-measurable on $[0,T]\times X$ and for a.e. $t\in(0,T)$ the mapping $v\mapsto
f(t,v)$ is lower semicontinuous on $X$. Note that if $f$ is a normal integrand, then by the  \textsc{Pettis} theorem, see, e.g., \textsc{Diestel \& Uhl} \cite[Theorem 2, p. 42]{DieUhl77VEME}, the mapping $t\mapsto f(t,v(t))$ is \textsc{Lebesgue} measurable for any \textsc{Bochner} measurable functional $v:[0,T]\rightarrow X$.

\begin{thm}\label{th:integral.functional} Let $X$ be a separable and reflexive \textsc{Banach} space and let $f:[0,T]\times X\rightarrow (-\infty,+\infty]$ be a normal integrand such that $f(t,\cdot): X\rightarrow (-\infty,+\infty]$ is a proper, lower semicontinuous and convex functional for a.e. $t\in (0,T)$ and denote with $F$ the integral functional defined by \eqref{ap:integral.functional}. Furthermore, with $f^*:[0,T]\times X^*\rightarrow (-\infty,+\infty]$ we denote the conjugate functional of $f$ given by $f^*(t,\cdot)= (f(t,\cdot))^*, t\in [0,T]$. Moreover, assume that there exist constants $\alpha,\alpha^*,\beta,\beta^*>0$ such that 
\begin{align*}
f(t,v)+\alpha \Vert v\Vert+\beta\geq 0 \quad \text{for a.e. }t\in [0,T] \text{ and all } v\in X,
\end{align*} and
\begin{align*}
f^*(t,\xi)+\alpha^* \Vert \xi\Vert_*+\beta^*\geq 0 \quad \text{for a.e. }t\in [0,T] \text{ and all } \xi\in X^*.
\end{align*} Then, the following assertions hold
\begin{itemize}
\item[$i)$] The functional $f^*:[0,T]\times X^*\rightarrow (-\infty,+\infty]$ is a normal integrand, and if $F$ is proper, then the conjugate functional $F^*:\rmL^{p^*}(0,T;X^*)\rightarrow \overline{\mathbb{R}}$ is proper, lower semicontinuous and convex, and is given by the integral functional
\begin{flalign*}
F^*(\xi)=\begin{cases}
\int_0^T f^*(t,\xi(t))\dd t \quad &\text{if } f^*(\cdot,\xi(\cdot))\in \rmL^1(0,T),\\
+\infty \quad &\text{otherwise}.
\end{cases}
\end{flalign*}
\item[$ii)$] The functional $F$ is lower semicontinuous and convex on $\rmL^p(0,T;X)$, and there holds $F(v)>-\infty$ for all $v\in \rmL^p(0,T;X)$.
\item[$iii)$] Let $F$ be proper, and let $v\in \DOM(F)$ and $\xi\in \rmL^{p^*}(0,T;X^*)$. Then, $\xi \in \partial F(v)\subset \rmL^{p^*}(0,T;X^*)$  if and only if $\xi(t)\in \partial f(t,v(t))\subset X^*\, \text{ for a.e. }t\in(0,T)$.
\end{itemize} 
\end{thm}

\begin{proof}[Proof]
Assertions $i)$ follows from \textsc{Kenmochi} \cite[Proposition 1.1]{Kenm75NPVI} and \textsc{Rockafellar}  \cite[Proposition 2 \& Theorem 2]{Rock71CIFD} and Assertion $ii)$ follows from \cite[Proposition 4.1 \& Corollary 4.1, p. 18]{EkeTem76CAVP}. Assertion $iii)$ follows from $i),ii)$, Lemma \ref{le:Leg.Fen}, and the fact that 
\begin{align}\label{eq:11}
\int_0^T \left( f(t,v(t))+f^*(t,\xi(t))-\langle \xi(t), v(t)\rangle_{X^*\times X}\right)\dd t=0 
\end{align} if and only if
\begin{align*}
f(t,v(t))+f^*(t,\xi(t))-\langle \xi(t), v(t)\rangle_{X^*\times X}=0\quad \text{a.e. in }(0,T),
\end{align*} which in turn follows from the fact that the integrand in \eqref{eq:11} is by the \textsc{Fenchel--Young} inequality, always non-negative.
\end{proof}

In the next result, we show the weak-weak closedness of the subdifferential of an integral function defined on a \textsc{Bochner} space. 
\begin{lem}\label{le:A1} Let the functionals $f, f_n: [0,T]\times X\rightarrow (-\infty,+\infty]$ be given and fulfill the assumptions of Theorem \ref{th:integral.functional}, and let $p\in (1,+\infty)$. Furthermore, let $(v_n)_{n\in \mathbb{N}}\subset \rmL^p(0,T;X)$ and $(\xi_n)_{n\in \mathbb{N}}\subset \rmL^{p^*}(0,T;X^*)$ with $\xi_n\in \partial F_n(v_n)$ such that $v_n\rightharpoonup v$ in $\rmL^p(0,T;X)$ and $\xi_n \rightharpoonup \xi$ in $\rmL^{p^*}(0,T;X^*)$ as $n\rightarrow \infty$, where $F_n$ is the integral functional associated to $f_n$. If 
\begin{align}\label{le:A1.assum.1}
 \int_{0}^T\left( f(t,v(t))+f^*(t,\xi(t))\right)\dd t\leq \liminf_{n\rightarrow \infty} \int_0^T\left( f_n(t,v_n(t))+f_n^*(t,\xi_n(t))\right)\dd t
\end{align} and 
\begin{align}\label{le:A1.assum.2}
\limsup_{n\rightarrow \infty}\int_0^T\langle \xi_n(t)-\xi(t), v_n(t)-v(t)\rangle_{X^*\times X}\dd t\leq 0,
\end{align} then $\xi(t)\in \partial f(t,v(t))$ a.e. in $(0,T)$ and there holds 
\begin{align*}
 \int_{0}^T\left( f(t,v(t))+f^*(t,\xi(t))\right)\dd t= \lim_{n\rightarrow \infty} \int_0^T\left( f_n(t,v_n(t))+f_n^*(t,\xi_n(t))\right)\dd t
\end{align*}
\end{lem}
\begin{proof}
By the \textsc{Legendre--Fenchel} inequality and Assumptions \ref{le:A1.assum.1} and \ref{le:A1.assum.2}, we find
\begin{align*}
 \int_{0}^T\left\langle \xi(t),v(t)\right\rangle_{X^*\times X}\dd t &\leq \int_{0}^T\left( f(t,v(t))+f^*(t,\xi(t))\right)\dd t\\
 &\leq \liminf_{n\rightarrow \infty} \int_0^T\left( f_n(t,v_n(t))+f_n^*(t,\xi_n(t))\right)\dd t\\
 &\leq \limsup_{n\rightarrow \infty} \int_0^T\left( f_n(t,v_n(t))+f_n^*(t,\xi_n(t))\right)\dd t\\
 &= \limsup_{n\rightarrow \infty} \int_0^T\left\langle \xi_n(t),v_n(t)\right\rangle_{X^*\times X}\dd t\\
 &= \int_{0}^T\left\langle \xi(t),v(t)\right \rangle_{X^*\times X}\dd t.
\end{align*} By Theorem \ref{th:integral.functional} and Lemma \ref{le:Leg.Fen}, it follows that $\xi(t)\in \partial f(t,v(t))\subset X^*\, \text{ for a.e. }t\in(0,T)$.
\end{proof}

\bibliographystyle{my_alpha}
\bibliography{alex_pub,bib_aras}

\end{document}